\documentclass[11pt]{article}
\usepackage{amsmath,amssymb,graphicx,epsfig,hyperref}

\newcommand{\qed}{\rule{3mm}{3mm}}

\newcommand{\zz}{{\frak z}}
\newcommand{\ee}{{\frak e}}

\newcommand{\cA}{{\cal A}}
\newcommand{\cG}{{\cal G}}
\newcommand{\cH}{{\cal H}}
\newcommand{\cL}{{\cal L}}

\newcommand{\cQ}{{\cal Q}}
\newcommand{\cT}{{\cal T}}

\newtheorem{theorem}{Theorem}[section]
\newtheorem{proposition}[theorem]{Proposition}
\newtheorem{lemma}[theorem]{Lemma}
\newtheorem{definition}[theorem]{Definition}
\newtheorem{corollary}[theorem]{Corollary}

\newcommand{\VTL}{V(\cT\cL)}
\newcommand{\VHL}{V(\cH\cL)}
\newcommand{\EHL}{E(\cH\cL)}
\newcommand{\FHL}{F(\cH\cL)}
\newcommand{\T}{T}
\newcommand{\R}{R}
\renewcommand{\S}{S}

\renewcommand{\bar}{\overline}
\newcommand{\A}{a}
\newcommand{\B}{b}
\newcommand{\C}{c}
\newenvironment{proof}{\noindent{\em Proof.} 
\normalfont}{\qed\par}
\DeclareRobustCommand{\qed}{%
  \ifmmode 
  \else \leavevmode\unskip\penalty9999 \hbox{}\nobreak\hfill
  \fi
  \quad\hbox{\qedsymbol}}
\newcommand{\openbox}{\leavevmode
  \hbox to.77778em{%
  \hfil\vrule
  \vbox to.675em{\hrule width.6em\vfil\hrule}%
  \vrule\hfil}}
\newcommand{\qedsymbol}{\openbox}
\newcommand{\quadmatrix}[4]{
  \left(
    \begin{array}{cc}
      #1&#2\\
      #3&#4
    \end{array}
  \right)}

\begin{document}

\begin{center}
{\large\bf
Hexagonal circle patterns and integrable systems.\\
Patterns with constant angles}
\end{center}
\bigskip

\begin{center}
{\sc Alexander I.\,Bobenko}\footnote{E--mail: {\tt bobenko}
{\makeatother @ \makeatletter}{\tt math.tu-berlin.de }} and {\sc
Tim \,Hoffmann}\footnote{E--mail: {\tt timh} {\makeatother @
\makeatletter}{\tt sfb288.math.tu-berlin.de }}
\end{center}
\begin{center}
Fakult\"at II, Institut f\"ur Mathematik, Technische Universit\"at Berlin, \\
Strasse des 17 Juni 136, 10623 Berlin, Germany
\end{center}
\bigskip




\section{Introduction}

The theory of circle packings and, more generally, of circle
patterns enjoys in recent years a fast development and a growing
interest of specialists in complex analysis and discrete
mathematics. This interest was initiated by Thurston's rediscovery
of the Koebe-Andreev theorem \cite{K} about circle packing
realizations of cell complexes of a prescribed combinatorics and
by his idea about approximating the Riemann mapping by circle
packings (see \cite{T1,RS}). Since then many other remarkable
facts about circle patterns were established, such as the discrete
maximum principle and Schwarz's lemma \cite{R} and the discrete
uniformization theorem \cite{BS}. These and other results
demonstrate surprisingly close analogy to the classical theory and
allow one to talk about an emerging of the "discrete analytic
function theory" \cite{DS}, containing the classical theory of
analytic functions as a small circles limit.

Approximation problems naturally lead to infinite circle patterns
for an analytic description of which it is advantageous to stick
with fixed regular combinatorics. The most popular are hexagonal
packings where each circle touches exactly six neighbors. The
$C^\infty$ convergence of these packings to the Riemann mapping
was established in \cite{HS}. Another interesting and elaborated
class with similar approximation properties to be mentioned here
are circle patterns with the combinatorics of the square grid
introduced by Schramm \cite{S}. The square grid combinatorics of
Schramm's patterns results in an analytic description which is
closer to the Cauchy-Riemann equations of complex analysis then
the one of the packings with hexagonal combinatorics. Various
other regular combinatorics also have similar properties \cite{H}.

Although computer experiments give convincing evidence for the
existence of circle packing analogs of many standard holomorphic
functions \cite{DS}, the only circle packings that have been
described explicitly are Doyle spirals \cite{BDS} (which are
analogs of the exponential function) and conformally symmetric
packings \cite{BH} (which are analogs of a quotient of Airy
functions). Schramm's patterns are richer with explicit examples:
discrete analogs of the functions $\exp (z), {\rm erf}(z)$, Airy
\cite{S} and $z^c, \log (z)$ \cite{AB} are known. Moreover ${\rm
erf}(z)$ is also an entire circle pattern.\footnote{Doyle
conjectured that the Doyle spirals are the only entire circle
packings. This conjecture remains open.}

A natural question is: what property is responsible for this
comparative richness of Schramm's patterns?  Is it due to the
packing - pattern or (hexagonal - square) combinatorics
difference? Or maybe it is the integrability of Schramm's patterns
which is crucial. Indeed, Schramm's square grid circle patterns in
conformal setting are known to be described by an integrable
system \cite{BP2} whereas for the packings it is still
unknown\footnote{It should be said that, generally, the subject of
discrete integrable systems on lattices different from ${\Bbb
Z}^n$ is underdeveloped at present. The list of relevant
publications is almost exhausted by \cite{Ad,KN,ND}.}.

In the present paper we introduce and study {\em hexagonal circle
patterns with constant angles}, which merge features of the two
circle patterns discussed above. Our circle patterns have the
combinatorics of the regular hexagonal lattice (i.e. of the
packings) and intersection properties of Schramm's patterns.
Moreover the latter are included as a special case into our class.
An example of a circle patterns studied in this paper is shown in
Figure \ref{fig:Doyle}. Each elementary hexagon of the honeycomb
lattice corresponds to a circle, and each common vertex of two
hexagons corresponds to an intersection point of the corresponding
circles. In particular, each circle carries six intersection
points with six neighboring circles and at each point there meet
three circles. To each of the three types of edges of the regular
hexagonal lattice (distinguished by their directions) we associate
an angle $0\leq\alpha_n<\pi,\ n=1,2,3$ and set that the
corresponding circles of the hexagonal circle pattern intersect at
this angle. It is easy to see that $\alpha_n$'s are subject to the
constraint $\alpha_1+\alpha_2+\alpha_3=\pi$.

We show that despite of the different combinatorics the properties
and the description of the hexagonal circle patterns with constant
angles are quite parallel to those of Schramm's circle patterns.
In particular, the intersection points of the circles are
described by a discrete equation of Toda type known to be
integrable \cite{Ad}. In Section \ref{s.Conformal} we present a
conformal (i.e. invariant with respect to M\"obius
transformations) description of the hexagonal circle patterns with
constant angles and show that one can vary the angles $\alpha_n$
arbitrarily preserving the cross-ratios of the intersection points
on circles, thus each circle pattern generates a two-parameter
deformation family. Analytic reformulation of this fact provides
us with a new integrable system possessing a Lax representation in
2 by 2 matrices on the regular hexagonal lattice.

 Conformally symmetric hexagonal
circle patterns are introduced in Section \ref{sec:conformal}.
Those are defined as patterns with conformally symmetric flowers,
i.e. each circle with its six neighbors is invariant under a
M\"obius involution (M\"obius $180^0$ rotation). A similar class
of circle packings was investigated in \cite{BH}. Let us mention
also that a different subclass of hexagonal circle patterns - with
the multi-ratio property instead of the angle condition - was
introduced and discussed in detail in \cite{BHS}. In particular it
was shown that this class is also described by an integrable
system. Conformally symmetric circle patterns comprise the
intersection set of the two known integrable classes of hexagonal
circle patterns: "with constant angles" and "with multi-ratio
property". The corresponding equations are linearizable and can be
easily solved.

Further in Section \ref{s.doyle} we establish a rather remarkable
fact. It turns out that Doyle circle packings and analogue
hexagonal circle patterns with constant angles are build out of
the same circles. Moreover, given such a pattern one can
arbitrarily vary the intersection angles $\alpha_n$ preserving the
circle radii.

Extending the intersection points of the circles by their centers
one embeds hexagonal circle patterns with constant angles into an
integrable system on the dual Kagome lattice (see Section
\ref{s.DualLax}). Having included hexagonal circle patterns with
constant angles into the framework of the theory of integrable
systems, we get an opportunity of applying the immense machinery
of the latter to study the properties of the former. This is
illustrated in Section \ref{s.z^c}, where we introduce and study
some isomonodromic solutions of our integrable systems on the dual
Kagome lattice. The corresponding circle patterns are natural
discrete versions of the analytic functions $z^c$ and $\log z$.
The results of Section \ref{s.z^c} constitute an extension to the
present, somewhat more intricate, situation of the similar
construction for Schramm's circle patterns with the combinatorics
of the square grid \cite{BP2,AB}.

\bigskip

{\bf Acknowledgements} The authors thank S.I.~Agafonov and
Yu.B.~Suris for collaboration and helpful discussions.

\section{Hexagonal circle patterns with constant angles}
\label{s.Hexagonal}

\begin{figure}[htbp]
  \begin{center}
\includegraphics[width=0.3\hsize]{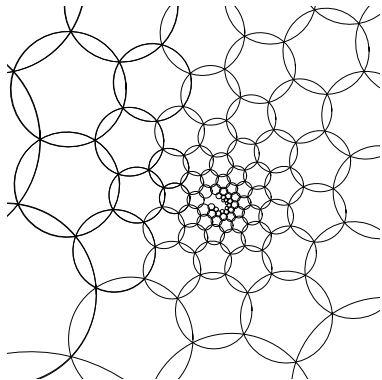}
    \caption{Doyle hexagonal isotropic circle pattern.}
    \label{fig:Doyle}
  \end{center}
\end{figure}

The present paper deals with {\em hexagonal circle patterns}, i.e.
circle patterns with the combinatorics of the regular hexagonal
lattice (the honeycomb lattice). An example is shown in Figure
\ref{fig:Doyle}. Each elementary hexagon of the honeycomb lattice
corresponds to a circle, and each common vertex of two hexagons
corresponds to an intersection point of the corresponding circles.
In particular, each circle carries six intersection points with
six neighboring circles and at each intersection point exactly
three circles meet.

For analytic description of hexagonal circle patterns we
introduce
some convenient lattices and variables. First of all we define the
{\em regular triangular lattice} $\cT\cL$ as the cell complex
whose vertices are
\begin{equation}
V(\cT\cL)=\Big\{\zz=k+\ell\omega+m\omega^2:\; k,\ell,m\in{\Bbb
Z}\Big\},\quad {\rm where}\quad \omega=\exp(2\pi i/3),
\end{equation}
whose edges are all non--ordered pairs
\begin{equation}
E(\cT\cL)=\Big\{[\zz_1,\zz_2]:\; \zz_1,\zz_2\in V(\cT\cL),\;
|\zz_1-\zz_2|=1 \Big\},
\end{equation}
and whose 2-cells are all regular triangles with the vertices in
$V(\cT\cL)$ and the edges in $E(\cT\cL)$. We shall use triples
$(k,\ell,m)\in {\Bbb Z}^3$ as coordinates of the vertices. On the
regular triangular lattice two such triples are equivalent and
should be identified if they differ by the vector $(n,n,n)$ with
$n\in{\Bbb Z}$.

The vertices of the regular triangular lattice correspond to
centers and intersection points of hexagonal circle patterns.
Associating one of the centers with the point $k=\ell=m=0$ we
obtain the {\em regular hexagonal sublattice} $\cH\cL$ with all
2-cells being the regular hexagons with the vertices in
\begin{equation}
V(\cH\cL)=\Big\{\zz=k+\ell\omega+m\omega^2:\; k,\ell,m\in{\Bbb
Z},\; k+\ell+m\not\equiv 0\!\!\pmod 3\Big\},
\end{equation}
and the edges in
\begin{equation}
E(\cH\cL)=\Big\{[\zz_1,\zz_2]:\;\zz_1,\zz_2\in V(\cH\cL),\;
|\zz_1-\zz_2|=1 \Big\}.
\end{equation}
The cells and the vertices correspond to circles and to the
intersection points of the hexagonal circle patterns respectively.
Natural labelling of the faces
$$
F(\cH\cL)=\Big\{\zz=k+\ell\omega+m\omega^2:\; k,\ell,m\in{\Bbb
Z},\; k+\ell+m\equiv 0\!\!\pmod 3\Big\}
$$
yields
$$ V(\cH\cL)  \bigcup F(\cH\cL) = V(\cT\cL).$$

\begin{definition}                              \label{def hex pattern}
We say that a map $w:V(\cH\cL)\mapsto\hat{\Bbb C}$ defines a {\em
hexagonal circle pattern}, if the following condition is
satisfied:
\begin{itemize}
\item Let
\[
\zz_k=\zz'+\varepsilon^k\in V(\cH\cL), \quad k=1,2,\ldots,6,\quad
where \quad\varepsilon=\exp(\pi i/3),
\]
 be the vertices of any elementary hexagon in $\cH\cL$ with the center
 $\zz'\in F(\cH\cL)$. Then the points
 $w(\zz_1),w(\zz_2),\ldots,w(\zz_6)\in\hat{\Bbb C}$ lie on a circle, and their
 circular order is just the listed one. We denote the circle through
 the points $w(\zz_1),w(\zz_2),\ldots,w(\zz_6)$ by $C(\zz')$, thus putting it
 into a correspondence with the center $\zz'$ of the elementary hexagon above.
 \end{itemize}
\end{definition}

As a consequence of this condition, we see that if two elementary
hexagons of $\cH\cL$ with the centers in $\zz',\zz''\in F(\cH\cL)$
have a common edge $[\zz_1,\zz_2]\in E(\cH\cL)$, then the circles
$C(\zz')$ and $C(\zz'')$ intersect in the points $w(\zz_1)$,
$w(\zz_2)$. Similarly, if three elementary hexagons of $\cH\cL$
with the centers in $\zz',\zz'',\zz''' \in F(\cH\cL)$ meet in one
point $\zz_0\in V(\cH\cL)$, then the circles $C(\zz')$, $C(\zz'')$
and $C(\zz''')$ also have a common intersection point $w(\zz_0)$.
\vspace{2mm}

{\bf Remark}. We will consider also circle patterns defined not on
the whole of $\cH\cL$, but rather on some connected subgraph of
the regular hexagonal lattice. \vspace{2mm}

To each pair of intersecting circles we associate the
corresponding edge $\ee\in E(\cH\cL)$ and denote by $\phi(\ee)$
the intersection angle of the circles, $0\le \phi< 2\pi$. The
edges of $E(\cH\cL)$ can be decomposed into three classes
\begin{eqnarray}                            \label{E^H123}
E^H_1&=&\Big\{\ee=[\zz',\zz'']\in E(\cH\cL):\,
\zz'-\zz''=\pm 1\Big\}, \nonumber\\
E^H_2&=&\Big\{\ee=[\zz',\zz'']\in E(\cH\cL):\,
\zz'-\zz''=\pm \omega\Big\}, \\
E^H_3&=&\Big\{\ee=[\zz',\zz'']\in E(\cH\cL):\,
\zz'-\zz''=\pm\omega^2\Big\}\nonumber.
\end{eqnarray}
These three sets correspond to three possible directions of the
edges of the regular hexagonal lattice.

We shall study in this paper a subclass of hexagonal circle
patterns intersecting at given angles defined globally on the
whole lattice.

\begin{definition}
We say that a map $w:V(\cH\cL)\mapsto\hat{\Bbb C}$ defines a {\em
hexagonal circle pattern with constant angles} if in addition to
the condition of Definition \ref{def hex pattern} the intersection
angles of circles are constant within their class $E^H_n$, i.e.
$$
\phi(\ee)=\alpha_n\qquad \forall \ee\in E^H_n,\, n=1,2,3.
$$
This angle condition implies in particular
\begin{equation}                    \label{sumalpha=pi}
\alpha_1+ \alpha_2+\alpha_3=\pi.
\end{equation}
\end{definition}

We call the circle pattern {\em isotropic} if all its intersection
angles are equal $\alpha_1=\alpha_2=\alpha_3=\pi/3$.

The existence of hexagonal circle patterns with constant angles
can be easily demonstrated via solving a suitable Cauchy problem.
For example, one can start with the initial circles
$C(n(1-\omega)), C(n(\omega-\omega^2)),
C(n(\omega^2-1)),\,n\in{\Bbb N}$.

{\bf Remark.} In the special case $\alpha_1=\alpha_2=\pi/2,
\alpha_3=0$ at each vertex one obtains two pairs of touching
circles intersecting orthogonally. The hexagonal circle pattern
becomes in this case a circle patterns of Schramm \cite{S} with
the combinatorics of the square grid. The generalization
$\alpha_1+\alpha_2=\pi, \alpha_3=0$ was introduced in \cite{BP2}.

\begin{figure}[htbp]
  \begin{center}
    \includegraphics[width=0.5\hsize]{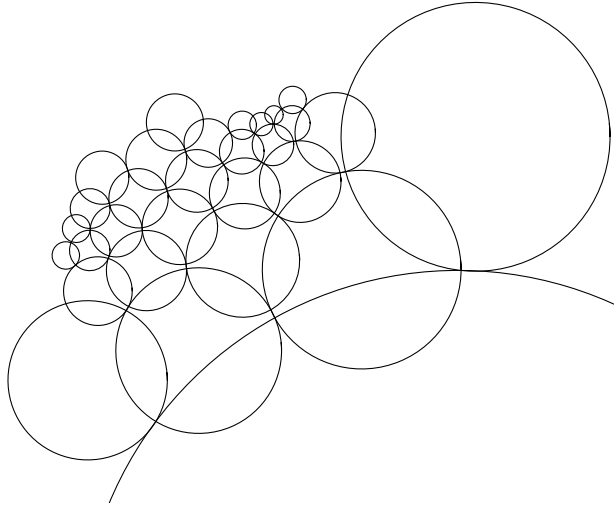}
    \caption{A nearly Schramm pattern.}
    \label{fig:NearlySchramm}
  \end{center}
\end{figure}

Figure~\ref{fig:NearlySchramm} shows a nearly Schramm pattern.
\vspace{2mm}

\section{Point and radii descriptions}\label{s.PointRadii}

In this paper three different analytic descriptions are used to
investigate hexagonal circle patterns with constant angles.
Obviously, these circle patterns can be characterized through the
{\em radii} of the circles. On the other hand, they can be
described through the coordinates of some natural {\em points},
such as the intersection points or centers of circles. Finally,
note that the class of circle patterns with constant angles is
invariant with respect to arbitrary fractional-linear
transformations of the Riemann sphere $\hat{\Bbb C}$ (M\"obius
transformations). Factorizing with respect to this group we
naturally come to a {\em conformal} description of the hexagonal
circle patterns with constant angles in Section \ref{s.Conformal}.

\begin{figure}[htbp]
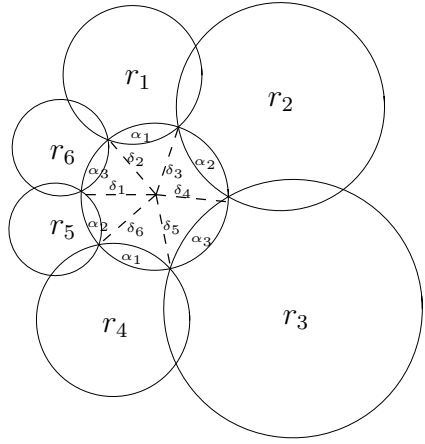

  \begin{center}
\input hexFlower.tex
    \caption{Circle flower.}
    \label{fig:flower}
  \end{center}
\end{figure}

The basic unit of a hexagonal circle pattern is the {\em flower},
illustrated in Figure \ref{fig:flower} and consisting of a center
circle surrounded by six petals. The radius $r$ of the central
circle and the radii $r_n,\, n=1,\ldots,6$ of the petals satisfy
\begin{equation}
\arg\prod_{n=1}^6 (r+e^{i\alpha_n}r_n) =\pi,
\end{equation}
where $\alpha_n$ is the angle between the circles with radii $r$
and $r_n$. Specifying this for the hexagonal circle patterns with
constant angles and using (\ref{sumalpha=pi}) one obtains the
following theorem.

\begin{theorem}
The mapping $r:F(\cH\cL) \to {\Bbb R}_+$ is the radius function of
a hexagonal circle pattern with constant angles $\alpha_1,
\alpha_2, \alpha_3$ iff it satisfies
\begin{equation}                        \label{eq.radii}
\arg\prod_{n=1}^3 (1+e^{i\alpha_n}R_n)(e^{-i\alpha_n}+R_{n+3}) =0,
\qquad R_n=\dfrac{r_n}{r}
\end{equation}
for every flower.
\end{theorem}
Conjugating (\ref{eq.radii}) and dividing it by the product
$R_1\dots R_6$ one observes that for every hexagonal circle
pattern with constant angles there exists a {\em dual} one.
\begin{definition}                      \label{d.dual}
Let $r:F(\cH\cL) \to {\Bbb R}_+$ be the radius function of a
hexagonal circle pattern $CP$ with constant angles. The hexagonal
circle pattern $CP^*$ with the same constant angles and the radii
function $r^*:F(\cH\cL) \to {\Bbb R}_+$ given by
\begin{equation}                        \label{dual.radii}
r^*=\dfrac{1}{r}
\end{equation}
is called {\em dual} to $CP$.
\end{definition}

For deriving the point description of the hexagonal circle
patterns with constant angles it is convenient to extend the
intersection points of the circles by their centers. Fix some
point $P_\infty\in\hat{\Bbb C}$. The reflections of this point in
the circles of the pattern will be called {\em conformal centers}
of the circles. In the particular case $P_\infty=\infty$ the
conformal centers become the centers of the corresponding circles.
We call an extension of a circle pattern by conformal centers a
{\em center extension}.

\begin{figure}[htbp]
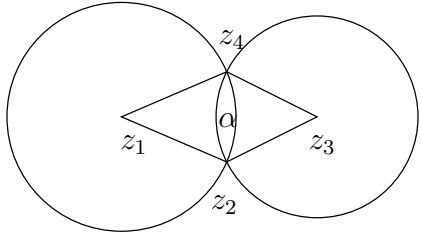

  \begin{center}
\input crossRatioAngle.tex
    \caption{Cross-ratio of an elementary quadrilateral.}
    \label{fig:cross-ratio}
  \end{center}
\end{figure}

M\"obius transformations play crucial role for the considerations
in this paper. Recall that  the {\em cross-ratio} of four points
\begin{equation}                            \label{cross-ratio}
q(z_1,z_2,z_3,z_4):=\dfrac{(z_2-z_1)(z_4-z_3)} {(z_3-z_2)(z_1-z_4)}
\end{equation}
is invariant with respect to these transformations. We start with
a simple
\begin{lemma}                               \label{l.cross-ratio=angle}
Let $z_2,z_4$ be the intersection points and $z_1,z_3$ the
conformal centers of two circles intersecting with the angle
$\alpha$ as in Figure \ref{fig:cross-ratio}. The cross-ratio of
these points is
\begin{equation}                            \label{cross-ratio=angle}
q(z_1,z_2,z_3,z_4)=e^{-2i\alpha}.
\end{equation}
\end{lemma}
The claim is obvious for the Euclidean centers. These can be
mapped to conformal centers by an appropriate M\"obius
transformation which preserves the cross-ratio.

Thus a circle pattern provides a solution $z:V(\cT\cL)\to\hat{\Bbb
C}$ to (\ref{cross-ratio=angle}) with the corresponding angles
$\alpha$. This solution is defined at the vertices of the lattice
$\cT\cL$ with quadrilateral sites comprised by the pairs of points
as in Lemma \ref{l.cross-ratio=angle}. For a hexagonal circle
pattern with constant angles one can derive from
(\ref{cross-ratio=angle}) equations for the intersection points
$z:V(\cH\cL)\to\hat{\Bbb C}$ and conformal centers
$z:V(\cT\cL\setminus\cH\cL)\to\hat{\Bbb C}$.

\begin{theorem}                                     \label{t.relToda}
Let $z, z_1, z_2, z_3, z_4, z_5, z_6$ be conformal centers of a
flower of a hexagonal circle pattern with constant angles
$\alpha_1, \alpha_2, \alpha_3$, where $\alpha_n, n=1,2,3$ are the
angles of pairs of circles corresponding to $z, z_n$ and $z,
z_{n+3}$. Define $\delta_1, \delta_2, \delta_3$ through
\begin{equation}                                    \label{eq.delta}
2\alpha_n=\delta_{n+2}-\delta_{n+1},\qquad \pmod {2\pi},\
n\in\{1,2,3\} \pmod 3.
\end{equation}
Then $z_n$ satisfy a discrete equation of Toda type on the
hexagonal lattice
\begin{equation}                                    \label{eq.relToda}
\sum_{n=1}^3 A_n
\left(\dfrac{1}{z-z_n}+\dfrac{1}{z-z_{n+3}}\right)=0,
\end{equation}
where
$$
A_n=e^{i\delta_{n+2}}- e^{i\delta_{n+1}} \qquad n\in\{1,2,3\}
\pmod 3.
$$
Let $w_1,w_2,w_3$ be the neighboring intersection points to the
point $w$ of a hexagonal circle pattern and let $\alpha_n$ be the
angle between the circles intersecting at $w,w_n$. Then the
following identity holds
\begin{equation}                                    \label{eq.w}
\sum_{n=1}^3 A_n \left(\dfrac{1}{w-w_n}\right)=0,
\end{equation}
\end{theorem}
\begin{figure}[htbp]
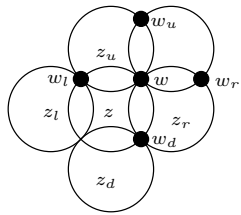

  \begin{center}
\input SchrammLabels.tex
    \caption{Intersection points and conformal centers of a Schramm circle pattern.}
    \label{fig:Schramm-labels}
  \end{center}
\end{figure}
In the special case of Schramm's patterns $\alpha_3=0,
\alpha_1+\alpha_2=\pi$ flowers contain only four petals and we
arrive at the following
\begin{theorem}                                 \label{t.Toda}
The intersection points $w, w_r, w_u, w_l, w_d\in {\Bbb C}$ and
the conformal centers $z, z_r, z_u, z_l, z_d\in {\Bbb C}$ of the
neighboring circles of a Schramm pattern (labelled as in Figure
\ref{fig:Schramm-labels}) satisfy the discrete equation of Toda
type on ${\Bbb Z}^2$ lattice
\begin{eqnarray}
\dfrac{1}{w-w_r}+ \dfrac{1}{w-w_l}&=&
\dfrac{1}{w-w_u}+\dfrac{1}{w-w_d},              \label{Toda_w}\\
\dfrac{1}{z-z_r}+ \dfrac{1}{z-z_{\ell}}&=&
\dfrac{1}{z-z_u}+\dfrac{1}{z-z_d}.              \label{Toda_z}
\end{eqnarray}
\end{theorem}


The proofs of these two theorems are presented in Appendix
\ref{Appendix.Toda} (in the case of complex $\alpha\in{\Bbb C}$).

It should be noticed that both equations (\ref{Toda_z}) and
(\ref{eq.relToda}) equations appeared in the theory of integrable
equations in a totally different context \cite{Su,Ad}. The
geometric interpretation in the present paper is new.

The sublattices of the centers and of the intersection points are
dependent and one can be essentially uniquely reconstructed from
the other. The corresponding formulas, which are natural to
generalize for complex cross-ratios (\ref{cross-ratio}), hold for
both discrete equations of Toda type discussed above. We present
these relations, which are of independent interest in the theory
of discrete integrable systems, in Appendix \ref{Appendix.Toda}.

In Section \ref{s.DualLax} we will show that the hexagonal circle
patterns with constant angles are described by an integrable
system on a regular lattice closely related to the lattices
introduced in Section \ref{s.Hexagonal}.

\section{Conformal description}         \label{s.Conformal}


Let us now turn our attention to a conformal description of the
hexagonal circle patterns. This description will be used in the
construction of conformally symmetric circle patterns in the next
section. We derive equations for the cross-ratios of the points of
the hexagonal lattice that allow to reconstruct the lattice up to
M\"obius transformations.

\begin{figure}[htbp]
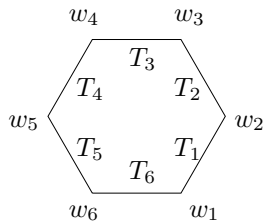

  \begin{center}
      \input Hexagon.tex
    \caption{Cross-ratios in a hexagon.}
    \label{fig:hexagon}
  \end{center}
\end{figure}

First we will investigate the relations of cross-ratios inside one
hexagon of the hexagonal lattice as shown in
Fig.~\ref{fig:hexagon}.
\begin{lemma}
  \label{thm:TLemma}
  Given a map $w:\VHL\to\hat{\mathbb C}$ let $\zz_1,\ldots,\zz_6$
  be six points of a hexagon cyclically ordered (see
  Fig.~\ref{fig:hexagon}).
  To each edge $[\zz_i,\zz_{i+1}]$ of the hexagon, let us assign a
  cross-ratio $\T_i$ of successive points $w_i = w(\zz_i)$:
  \begin{equation}
    \label{eq:Tdefinition}
    \T_i := q(w_i,w_{i+1},w_{i+2},w_{i-1}), \qquad (i\!\!\!\mod 6).
  \end{equation}
  Then the equations
  \begin{equation}
    \label{eq:Tequation}
    \frac{\T_1}{\T_4} = \frac{\T_3}{\T_6} = \frac{\T_5}{\T_2} = \frac{ \T_1 +
    \T_3 -1 - \T_1 T_2 \T_3}{ 1 - \T_2}
  \end{equation}
  hold.
\end{lemma}
\begin{proof}
  Let $m_i$ be the M\"obius transformation that maps
  $w_{i-1}, w_i$, and $w_{i+1}$ to 0, 1, and $\infty$. Then $M_i :=
  m_{i+1}^{-1}  m_i$ maps $w_{i-1}, w_i$, and $w_{i+1}$ to $w_{i},
  w_{i+1}$, and $w_{i+2 }$ respectively and has the form
  \begin{equation*}
    M_i = \left(
      \begin{array}{cc}
        1 &   -1\\
        1 & -T_i
      \end{array}
    \right)
  \end{equation*}
  Now $(M_6 M_5 M_4)^{-1} = \rho M_3 M_2 M_1$ for some $\rho$ gives
  the desired identities.
\end{proof}

Since every edge in $\EHL$ belongs to two hexagons it carries two
cross-ratios in general. We will now investigate the relation
between them and show that in case of a hexagonal circle pattern
with constant angles the two cross-ratios coincide. This way
Lemma~\ref{thm:TLemma} will furnish a map $\T:\EHL\to\mathbb C$.

Besides the cross-ratios of successive points in each hexagon (the
$\T_i$) we will need the cross-ratios of a point and its three
neighbors:
\begin{definition}
  Let $\zz_1,\zz_{2},$ and $\zz_{3}$ be the neighbors of $\zz\in\VHL$
  counterclockwise ordered and $[\zz_i,\zz]\in E^H_i, i = 1,2,3$. Any map
  $w:\VHL\to\hat{\mathbb C}$ furnishes three cross-ratios in each
  point $\zz\in \VHL$:
  \begin{equation}
    \label{eq:defS}
    \S_{\zz}^{(i)} :=
    q(w_i,w_{i+1},w,w_{i+2}),
    \quad i = 1,2,3 \quad (\!\!\!\mod 3) .
  \end{equation}
  They are linked by the modular transformation
\begin{equation}
  \label{eq:SRelation}
  \S^{(i+1)}_{\zz} = 1 - \frac1{\S^{(i)}_{\zz}}, \qquad (i\!\!\!\mod 3)
\end{equation}
\end{definition}

\begin{figure}[htbp]
  \begin{center}
    \input{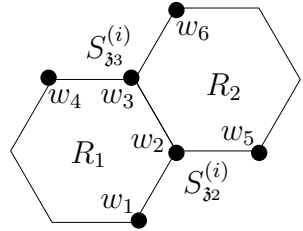}
    \caption{Six points around one edge.}
    \label{fig:TwoHexagons}
  \end{center}
\end{figure}
Of course the two types of cross-ratios are not independent: When
we look at an edge from $\EHL$ with its four neighboring edges the
six points form four cross-ratios (two $\T$ and two $\S$). Three
points fix the M\"obius transformation and the other three can be
calculated from only three cross-ratios. Therefore one expects one
equation:
\begin{lemma}
  Given $w:\VHL\to\hat{\mathbb C}$ and $\zz_1,\ldots,\zz_6$ as shown in
  Fig.~\ref{fig:TwoHexagons} and the cross-ratios
  \begin{equation}
    \R_1 := q(w_1,w_2,w_4,w_3), \quad \R_2 := q(w_5,w_2,w_6,w_3)
  \end{equation}
  \begin{equation}
    \label{eq:defS2}
    \S_{\zz_3}^{(i)} := q(w_3,w_4,w_2,w_6),\quad
    \S_{\zz_2}^{(i)} := q(w_2,w_5,w_3,w_1)
  \end{equation}
  with $[\zz_2,\zz_3]\in E^H_i$. Then the following identity holds:
  \begin{equation}
    \label{eq:leftRightRelation}
    \frac{1-\T_2}{1-\T_1} = \frac{\R_1}{\R_2} =
    \frac{\S_{\zz_3}^{(i)}}{\S_{\zz_2}^{(i)}}.
  \end{equation}
\end{lemma}
\begin{proof}
Use that $\T_i$ and $\R_i$ are linked by the relation
\begin{equation}
  \label{eq:defR}
  \R_i := \frac1{1-\T_i} = q(w_i,w_{i+1},w_{i-1},w_{i+2})
\end{equation}
and insert the definition of the cross-ratio.
\end{proof}

Everything written so far holds for any hexagonal lattice, but in
case of the hexagonal circle patterns with constant angles we can
calculate all $\S$ in terms of the $\alpha_i$ only:
\begin{equation}
  \label{eq:circleS}
  \S_\zz^{(i)} := q(w_i,w_{i+1},w,w_{i+2}) = e^{-i\alpha_{i}}\frac{\sin
  \alpha_{i+1}}{\sin \alpha_{i+2}}
\end{equation}
if $[\zz_i,\zz]\in E_i$. This can be verified easily by applying a
M\"obius transformation sending $\zz$ to infinity. The three
circles become straight lines forming a triangle with angles $\pi
- \alpha_i$. The $\S^{(i)}$ are now the quotients of two of its
edges.

In particular we get $\S^{(i)}_{\zz_1} = \S^{(i)}_{\zz_2}$ along
all edges $[\zz_1,\zz_2]\in\EHL$. Thus for our hexagonal circle
patterns the above formula~(\ref{eq:circleS}) implies that there
is only one $\T$ per edge since we get from equation
(\ref{eq:leftRightRelation})
\begin{equation}
  \label{eq:circleT}
  \T_1 = \T_2.
\end{equation}
So Lemma~\ref{thm:TLemma} defines in fact a map $\T:\EHL\to
\hat{\mathbb C}$.

Note that we can get back the $\alpha_i$ from say $S^{(1)}$ in
equation (\ref{eq:circleS}) by the following formulas:
\begin{equation}
  \label{eq:anglesFromS}
  e^{2 i\alpha_1} =\frac{\bar S^{(1)}}{S^{(1)}},\quad
  e^{2 i\alpha_2} = \frac{1+S^{(1)}}{1 +\bar S^{(1)}},\quad e^{2
  i\alpha_3} = \frac{S^{(1)}}{\bar S^{(1)}}\; \frac{1 +\bar S^{(1)}}{1+S^{(1)}}
\end{equation}

Now we can formulate the main theorem that states that the
equations which link the cross-ratios in a hexagon plus the
constantness of $\S^{(1)}$ on the hexagonal lattice describe for
real-negative $\T$ a hexagonal circle pattern with constant angles
up to M\"obius transformations:
\begin{theorem}
\label{thm:mainConformalThm} \noindent
  \begin{enumerate}
  \item
    Given $\S\in\mathbb C$ and a map $w:\VHL\to\hat{\mathbb C}$ for
    which $\S^{(1)}_{\zz}= \S$ for all $\zz\in\VHL$.
    Then the cross-ratios $T_i$
    of Lemma~\ref{thm:TLemma} define a map $\T:\EHL\to\hat{\mathbb C}$
    and obey equations (\ref{eq:Tequation}).
  \item
    Conversely given a solution $\T:\EHL\to\hat{\mathbb C}$ of
    (\ref{eq:Tequation}) then for each $\S\in{\mathbb C}$ there
    is up to M\"obius transformations a unique map
    $w:\VHL\to\hat{\mathbb C}$ having the map $\T$ as cross-ratios in the
    sense of Lemma~\ref{thm:TLemma} and having $\S^{(1)}_{\zz} = \S$ for
    all $\zz\in\VHL$.

    If the solution $\T$ is real-negative the resulting $w$ define a
    hexagonal circle pattern with constant angles.
  \end{enumerate}
\end{theorem}
\begin{proof}
  1. is proven above.\\
  2. Starting with the points $w_{2,0,0}, w_{1,1,0}$, and $w_{1,0,1}$
  (which fixes the M\"obius transformation) we can calculate
  $w_{1,0,0}$ using $\S^{(1)}_{1,0,0}$. Now we have three points for
  each hexagon touching in $\zz_{1,0,0}$. Using the $\T$'s we can
  determine all other points of them. Every new point has at least two
  determined neighbors so we can use the $\S^{(1)}$ to compute the
  third. Now we are able to compute all points of the new neighboring
  hexagons and so on.

  Since real-negative cross-ratios imply that the four points lie
  cyclically ordered on a circle, we have a circle pattern in the case
  of real-negative map $\T$. But with (\ref{eq:circleS}) the
  constantness of $S^{(1)}$ implies that the intersection angles are
  constant too.
\end{proof}

We have shown, that hexagonal circle patterns with constant angles
come in one complex (or two real) parameter families: one can
choose $\S\in\mathbb C$ arbitrarily preserving $\T$'s. In
Appendix~\ref{sec:confLax} it is shown how this observation
implies a Lax representation on the hexagonal lattice for the
system (\ref{eq:Tequation}).

Figure~\ref{fig:NearlySchramm} shows a nearly Schramm pattern. One
can obtain his description by taking combinations of $\T_i$'s and
$\S^{(i)}$'s that stay finite in the limit $\alpha_3\to 0$.

\section{Conformally symmetric circle patterns}\label{sec:conformal}

The basic notion of conformal symmetry introduced in \cite{BH} for
circle packings can be easily generalized to circle patterns:
Every elementary flower shall be invariant under the M\"obius
equivalent of a $180^\circ$ rotation.

\begin{definition}
  \begin{enumerate}
  \item An elementary flower of a hexagonal circle pattern with petals
    $C_i$ is called {\em conformally symmetric} if there exists a
    M\"obius involution sending $C_i$ to $C_{i+3}$ $(i\!\!\!\mod 3)$.
  \item A hexagonal circle pattern is called {\em conformally
      symmetric} if all of its elementary flowers are.
  \end{enumerate}
\end{definition}

For investigation of conformally symmetric patterns we will need
the notion of the multi-ratio of six points:
\begin{equation}
  \label{eq:muti-ratio}
  m(z_1,z_2,z_3,z_4,z_5,z_6) := \frac{z_1-z_2}{z_2-z_3}
  \frac{z_3-z_4}{z_4-z_5} \frac{z_5-z_6}{z_6-z_1}.
\end{equation}
Hexagonal circle pattern with multi-ratio $-1$ are discussed in
\cite{BHS} (see \cite{KS} for further geometric interpretation of
this quantity). It turns out that in the case of conformally
symmetric hexagonal circle patterns the two known integrable
classes ``with constant angles'' and ``with multi-ratio $-1$''
coincide.


\begin{proposition}
  \begin{enumerate}
  \item An elementary flower of a hexagonal circle pattern is
    conformally symmetric if and only if the opposite
    intersection angles are equal and the six intersection
    points with the central circle have multi-ratio $-1$.
  \item A hexagonal circle pattern is conformally symmetric
    iff and only if it has constant intersection angles and for all
    circles the six intersection points have multi-ratio $-1$.
  \end{enumerate}
\end{proposition}
\begin{proof}
  First let us show that six points $z_i$ have multi-ratio $-1$ if
  and only if there is a M\"obuis involution sending $z_i$ to
  $z_{i+3}$.

  If there is such M\"obius transformation it is clear, that
  $q(z_1,z_2,z_3,z_4) = q(z_4,z_5,z_6,z_1)$ and
  \begin{equation}
    \label{eq:mrCr}
    m(z_1,z_2,z_3,z_4,z_5,z_6) =
    -q(z_1,z_2,z_3,z_4)/q(z_4,z_5,z_6,z_1).
  \end{equation}
  implies
  \begin{equation}
    \label{eq:mrmr}
m(z_1,z_2,z_3,z_4,z_5,z_6) = -1.
  \end{equation}
  Conversely let $M$ be the M\"obius transformation sending $z_1, z_2$
  and $z_3$ to $z_4,z_5$ and $z_6$. Then for $z_* := M(z_4)$
  equation~(\ref{eq:mrmr}) implies $q(z_4,z_5,z_6,z_1) =
  q(z_1,z_2,z_3,z_4) = q(M(z_1),M(z_2),M(z_3),M(z_4)) =
  q(z_4,z_5,z_6,z_*)$ and thus $\zz_*=\zz_1$. The same computation
  yields $\zz_2 = M(\zz_5)$ and $\zz_3 = M(\zz_6)$.

  Now the first statement of the theorem is proven, since the
  intersection points and angles determine the petals completely. For
  the prof of the second statement the only thing left to show is that
  all flowers being conformally symmetric implies that the three
  intersection angles per flower sum up to $\pi$. So let us look at a
  flower around the circle $C$ with petals $C_i$. Let the angle
  between $C$ and $C_i$ be $\alpha_i$ (and we know that $\alpha_i =
  \alpha_{i+3}$). Then the angle $\beta_1 = \pi-\alpha_1-\alpha_2$ is
  the angle between $C_1$ and $C_2$ and $\beta_3 =
  \pi-\alpha_3-\alpha_1$ the one between $C_3$ and $C_4$. Since the
  flowers around $C_2$ and $C_3$ are conformally symmetric too we now
  have two ways to compute the angle $\beta_2$ between $C_2$ and
  $C_3$. Namely
  \begin{equation*}
    \pi-\alpha_2 -\alpha_3 = \beta_2 = \pi -(\pi-\alpha_3-\alpha_1)
    -(\pi-\alpha_1-\alpha_2)
  \end{equation*}
  which implies $\alpha_1+\alpha_2+\alpha_3 = \pi$.
\end{proof}

Using (\ref{eq:mrCr}) we see, that in the case of multi-ratio $-1$
the opposite cross-ratios $T$ defined in section~\ref{s.Conformal}
must be equal: $\T_i = \T_{(i\!\!\!\mod3)}$. Thus on $\EHL$ the
$T$'s must be constant in the direction perpendicular to the edge
they are associated with.  For the three cross-ratios in a hexagon
we get from (\ref{eq:Tequation})
\begin{equation}
  \label{eq:TequationSymmetric}
  \T_1 + \T_2 + \T_3 - \T_1 \T_2 \T_3 = 2.
\end{equation}
\begin{figure}[htbp]
  \begin{center}
    \input{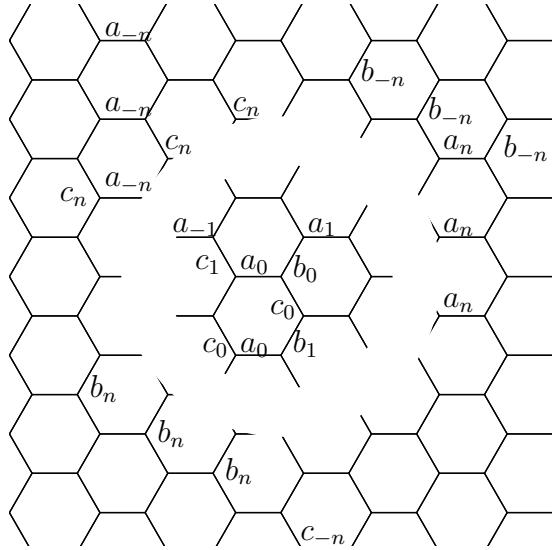}
    \caption{Cross-ratios for conformally symmetric circle patterns}
    \label{fig:gridLabels}
  \end{center}
\end{figure}
Let us rewrite this equation by using the labelling shown in
Fig.~\ref{fig:gridLabels}.
We denote by $\A_k, \B_\ell$ and $\C_m$ the cross-ratios
(\ref{eq:defR}) associated with the edges of the families
$E^H_1,E^H_2$ and $E^H_3$ respectively. Note that for the labels
in Fig.~\ref{fig:gridLabels}
$$ k +\ell+m = 1$$
holds. Written in terms of $\A_k,\B_\ell,\C_m$
equation~(\ref{eq:TequationSymmetric}) is linear
\begin{equation}
  \label{eq:ABCequation}
  \A_k + \B_\ell + \C_m = 1
\end{equation}
and can be solved explicitly:

\begin{lemma}
\label{thm:ABClemma}
  The general solution to (\ref{eq:ABCequation}) on $\EHL$ is given by
  \begin{equation}
    \label{eq:ABCsolution}
    \begin{array}{rcl}
      \A_k &=& \A_0 + k \Delta\\
      \B_\ell &=& \B_0 + \ell \Delta\\
      \C_m &=& \C_0 + m \Delta
    \end{array}
  \end{equation}
  with  $\A_0,\B_0,\C_0\in {\mathbb C}$ and $\Delta = 1 -\A_0
  -\B_0-\C_0$.
\end{lemma}
\begin{proof}
  Obviously (\ref{eq:ABCsolution}) solves (\ref{eq:ABCequation}). On
  the other hand  it is easy to
  show, that given the cross-ratios on three neighboring edges
  determine all other cross-ratios recursively. Therefore
  (\ref{eq:ABCsolution}) is the only solution to (\ref{eq:ABCequation})
\end{proof}
\begin{figure}[htbp]
  \begin{center}
    \includegraphics[width=0.5\hsize]{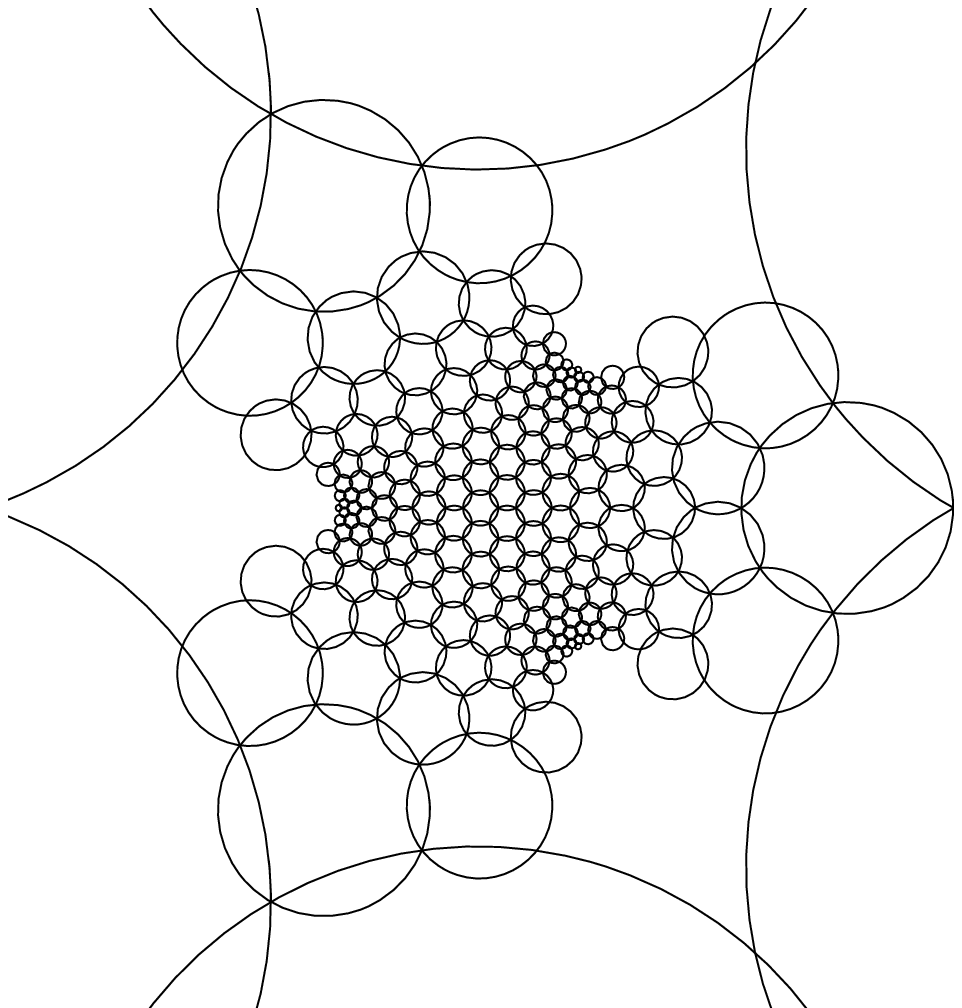}\includegraphics[width=0.5\hsize]{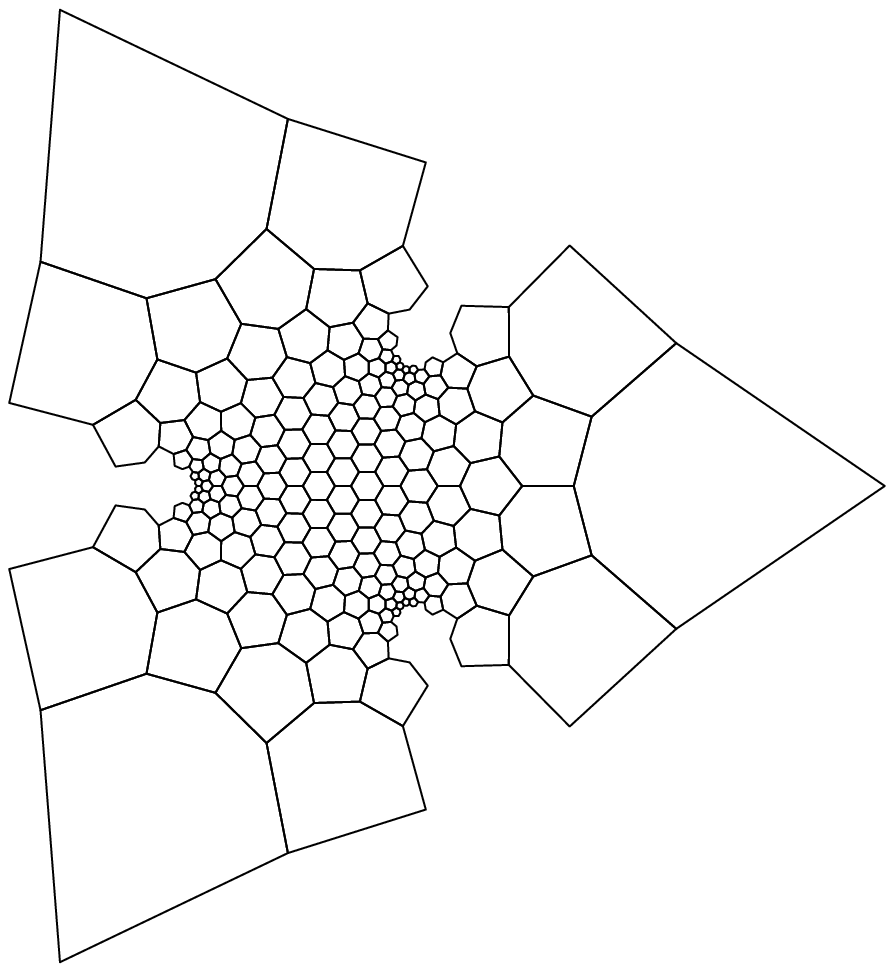}
    \caption{A conformally symmetric circle pattern and $\VHL$ under
      the quotient of two Airy functions $f(z) = \frac{{\rm Bi}(z) +
        \sqrt{3} {\rm Ai}(z)}{{\rm Bi}(z) - \sqrt{3} {\rm
        Ai}(z)}$.}
    \label{fig:Airy}
  \end{center}
\end{figure}
\begin{theorem}
  Conformally symetric circle hexagonal circle patterns are described
  as follows:
  Given $\A_n,\B_n,\C_n$ by (\ref{eq:ABCsolution}), choose
  $\S^{(1)}\in{\mathbb C}$ (or angles $\alpha_1$ and $\alpha_2$) then there
  is a conformally symmetric circle pattern with intersection angles
  $\alpha_1, \alpha_2$, and $\alpha_3$ given by formula (\ref{eq:anglesFromS}).
\end{theorem}

The cross-ratio of four points can be viewed as a discretization
of the Schwarzian derivative. In this sense conformally symmetric
patterns correspond to maps with a linear Schwarzian. The latter
are quotients of two Airy functions \cite{BH}.
Figure~\ref{fig:Airy} shows that a symmetric solution of
(\ref{eq:ABCequation}) is a good approximation of its smooth
counterpart. Figure~\ref{fig:NearlySchramm} also shows a
conformally symmetric pattern with $\S^{(1)} = 10i$ and $\A_0 =
\B_0 = 1/3 - 0.29$, and $\C_0 = 1-2 \A_0$.

\section{Doyle circle patterns}
\label{s.doyle} Doyle circle packings are described through their
radii function. The elementary flower of a hexagonal circle
packing is a central circle with six touching petals (which in
turn touch each other cyclically). Let the radius of the central
circle be $R$ and $R_1,\ldots,R_6$ the radii of the petals. Doyle
spirals are described through the constraint
\begin{equation}
  \label{eq:doyleRadii}
  R_k R_{k+3} = R^2, \quad R_k R_{k+2} R_{k+4} = R^3.
\end{equation}
There are two degrees of freedom for the whole packing---e.~g.\
$R_1/R$ and $R_2/R$ which are constant for all flowers. The next
lemma claims that the same circles that form a Doyle packing build
up a circle pattern with constant angles.
\begin{lemma}
  The radii function of a Doyle packing (i.~e.\ a solution to
  (\ref{eq:doyleRadii})) solves the pattern radii equation
  (\ref{eq.radii}) for any choice of the angles
  $\alpha_1,\alpha_2, \alpha_3 = \pi-\alpha_1-\alpha_2$.
\end{lemma}

\begin{proof}
  Now insert identities~(\ref{eq:doyleRadii}) into (\ref{eq.radii}).
\end{proof}
\begin{definition}
  A hexagonal circle pattern with constant angles who's radii function
  obeys the constraint (\ref{eq:doyleRadii}) is called a Doyle
  pattern.
\end{definition}
Figure~\ref{fig:Doyle} shows a Doyle pattern. The following lemma
and theorem shows how the Doyle patterns fit into our conformal
description:
\begin{lemma}
  Doyle patterns are conformally symmetric.
\end{lemma}
\begin{proof}
  We have to show that for Doyle pattern the multi-ratio of each
  hexagon is $-1$. One can assume that the circumfencing circle has
  radius 1 and center 0 and that $w_1 =1$. Then the other points are
  given by
  $$w_j = w_{j-1} \frac{1 + R_j e^{i\alpha_{(j mod 3)}}}%
  {1 + R_j e^{-i\alpha_{(j mod 3)}}}$$
  where $R_j$ are the radii of the petals. Inserting this into the
  definition of the multi-ratio implies the claim.
\end{proof}
\begin{theorem}
  Doyle patterns and their M\"obius transforms can be characterized in
  the following way: The corresponding solution to
  (\ref{eq:ABCequation}) is constant\footnote{Constant means here
    $\A_k=\A_0, \B_k=\B_0, \C_k=\C_0$ for all $k\in{\mathbb Z}$.},
  i.~e.\ $\A_0 +\B_0+\C_0 = 1$ and $\A_0,\B_0,\C_0>0$.
\end{theorem}
\begin{proof}
  Any Doyle pattern gives rise to a constant solution of
  (\ref{eq:ABCequation}) since all elementary flowers of a Doyle
  pattern are similar.

  On the other hand one sees easily, that all $0<\R_i<1$ can be realized.

\end{proof}

\section{Lax representation and Dual patterns}\label{s.DualLax}

We start with a general construction of ``integrable systems'' on
graphs which does not hang on the specific features of the
lattice. This notion includes the following ingredients:
\begin{itemize}
\item An {\it oriented graph} $\cG$ with the vertices $V(\cG)$ and
the edges $E(\cG)$.
\item A loop group $G[\lambda]$, whose elements are functions from
${\Bbb C}$ into some group $G$. The complex argument $\lambda$ of
these functions is known in the theory of integrable systems as
the ``spectral parameter''.
\item A ``wave function'' $\Psi: V(\cG)\mapsto G[\lambda]$,
defined on the vertices of $\cG$.
\item A collection of ``transition matrices'' $L: E(\cG)\mapsto G[\lambda]$
defined on the edges of $\cG$.
\end{itemize}
It is supposed that for any oriented edge
$\ee=(\zz_{out},\zz_{in})\in E(\cG)$ the values of the wave
functions in its ends are connected via
\begin{equation}                        \label{wave function evol}
\Psi(\zz_{in},\lambda)=L(\ee,\lambda)\Psi(\zz_{out},\lambda).
\end{equation}
Therefore the following {\it zero curvature condition} has to be
satisfied. Consider any closed contour consisting of a finite
number of edges of $\cG$:
\[
\ee_1=(\zz_1,\zz_2),\quad \ee_2=(\zz_2,\zz_3),\quad \ldots,\quad
\ee_p=(\zz_p,\zz_1).
\]
Then
\begin{equation}                            \label{zero curv cond}
L(\ee_p,\lambda)\cdots L(\ee_2,\lambda)L(\ee_1,\lambda)=I.
\end{equation}
In particular, for any edge $\ee=(\zz_1,\zz_2)$ one has
$\ee^{-1}=(\zz_2,\zz_1)$ and
\begin{equation}\label{zero curv cond inv}
L(\ee^{-1},\lambda)=\Big(L(\ee,\lambda)\Big)^{-1}.
\end{equation}

Actually, in applications the matrices $L(\ee,\lambda)$ depend
also on a point of some set $X$ (the ``phase space'' of an
integrable system), so that some elements $x(\ee)\in X$ are
attached to the edges $\ee$ of $\cG$. In this case the discrete
zero curvature condition (\ref{zero curv cond}) becomes equivalent
to the collection of equations relating the fields $x(\ee_1)$,
$\ldots$, $x(\ee_p)$ attached to the edges of each closed contour.
We say that this collection of equations admits a {\it zero
curvature representation}. Such representation may be used to
apply analytic methods for finding concrete solutions,
transformations or conserved quantities. \vspace{3mm}

\begin{figure}[htbp]
  \begin{center}
\input{quadLattice.tex}
    \caption{Quadrilateral lattice $\cQ\cL$.}
    \label{fig:QLattice}
  \end{center}
\end{figure}
In this paper we will deal with zero curvature representations on
the hexagonal lattice $\cH\cL$ and especially on a closely related
to it a {\em special quadrilateral lattice} $\cQ\cL$ which is
obtained from $\cH\cL$ by deleting from $E(\cT\cL)$ the edges of
the hexagonal lattice $E(\cH\cL)$ (i.e. corresponding to the
intersection points). So defined lattice $\cQ\cL$ has
quadrilateral cells, vertices $V(\cQ\cL)=V(\cT\cL)$ and the edges
\begin{equation*}
E(\cQ\cL)=E(\cT\cL)\setminus E(\cH\cL)
\end{equation*}
as shown in Figure \ref{fig:QLattice}. The lattice dual to
$\cQ\cL$ is known as Kagome lattice \cite{B}. Similar to
(\ref{E^H123}) there are tree types of edges in $E(\cQ\cL)$
distinguished by their directions:
\begin{eqnarray}                            \label{E^Q123}
E^Q_1&=&\Big\{\ee=[\zz',\zz'']\in E(\cQ\cL):\,
\zz'-\zz''=\pm 1\Big\}, \nonumber\\
E^Q_2&=&\Big\{\ee=[\zz',\zz'']\in E(\cQ\cL):\,
\zz'-\zz''=\pm \omega\Big\}, \\
E^Q_3&=&\Big\{\ee=[\zz',\zz'']\in E(\cQ\cL):\, \zz'-\zz''=
\pm\omega^2\Big\}\nonumber.
\end{eqnarray}

There is a natural labelling $(k,\ell,m)\in {\Bbb Z}^3$ of the
vertices $V(\cQ\cL)$ which respects the lattice structure of
$\cQ\cL$. Let us decompose $V(\cQ\cL)=V_0\cup V_1\cup V_{-1}$ into
three sublattices
\begin{eqnarray}                        \label{klm normalization}
V_0&=&\Big\{\zz=k+\ell\omega+m\omega^2:\; k,\ell,m\in{\Bbb Z},
\ k+\ell+m=0 \Big\},\nonumber\\
V_1&=&\Big\{\zz=k+\ell\omega+m\omega^2:\; k,\ell,m\in{\Bbb Z},
\ k+\ell+m=1 \Big\},\\
V_{-1}&=&\Big\{\zz=k+\ell\omega+m\omega^2:\; k,\ell,m\in{\Bbb Z},
\ k+\ell+m=-1 \Big\}\nonumber.
\end{eqnarray}
Note that this definition associates a unique triple $(k,\ell,m)$
to each vertex. Neighboring vertices of $\cQ\cL$, i.e. those
connected by edges, are characterized by the property that their
$(k,\ell,m)$-labels differ only in one component. Each vertex of
$V_0$ has six edges whereas the vertices of $V_{\pm 1}$ have only
three.

The $(k,\ell,m)$ labelling of the vertices of $\cQ\cL$ suggests to
consider this lattice as the intersection of ${\Bbb Z}^3$ with the
strip $\mid k+\ell+m \mid \leq 1$. This description turns out to
be useful especially for construction of discrete analogues of
$z^c$ and $\log z$ in Section \ref{s.z^c}. We denote
$p=(k,\ell,m)\in{\Bbb Z}^3$ and three different types of edges of
${\Bbb Z}^3$ by
\begin{eqnarray*}
E_1=(p,p+(1,0,0)),\ E_2=(p,p+(0,1,0)),\ E_3=(p,p+(0,0,1)),\
p\in{\Bbb Z}^3.
\end{eqnarray*}

The group $G[\lambda]$ we use in our construction is the {\it
twisted loop group} over ${\rm SL}(2,{\Bbb C})$:
\begin{equation}\label{loop group}
G[\lambda]=\Big\{L:{\Bbb C}\mapsto {\rm SL}(2,{\Bbb C})\Big|\;
L(-\lambda)=\sigma_3 L(\lambda)\sigma_3\Big\},\qquad \sigma_3={\rm
diag}(1,-1).
\end{equation}
To each type $E_n$ we associate a constant $\Delta_n\in{\Bbb C}$.
Let $z_{k,\ell,m},z^*_{k,\ell,m}$ be fields defined at vertices
$z,z^*:{\Bbb Z}^3\to {\Bbb C}$ and satisfying
\begin{equation}                            \label{dual}
(z_{in}-z_{out})(z^*_{in}-z^*_{out})=\Delta_n
\end{equation}
on all the edges $\ee=(p_{out},p_{in})\in E_n, n=1,2,3$. Here we
denote $z_{out, in}=z(p_{out, in}), z^*_{out, in}=z^*(p_{out,
in})$. We call the mapping $z^*$ {\em dual} to $z$. To each
oriented edge $\ee=(p_{out},p_{in})\in E_n$ we attach the
following element of the group $G[\lambda]$:
\begin{equation}                            \label{L}
L^{(n)}(\lambda)=(1-\lambda^2
\Delta_n)^{-1/2}\left(\begin{array}{cc}
1 & \lambda (z_{in}-z_{out})\\
\lambda (z^*_{in}-z^*_{out}) & 1 \end{array}\right).
\end{equation}
Note that this form as well as the condition (\ref{dual}) are
independent of the orientation of the edge. Substituting
(\ref{dual}) into (\ref{L}) one obtains $L(\lambda)$ in terms of
the field $z$ only.

Dual fields $z,z^*$ can be characterized in their own terms. The
condition that for any quadrilateral of the dual lattice the
oriented edges $(z^*_{out},z^*_{in})$ defined by (\ref{dual}) sum
up to zero imply the following statements.
\begin{theorem}                     \label{t.cross-ratio delta/delta}
(i) Let $z,z^*:{\Bbb Z}^3\to{\Bbb C}$ be dual fields, i.e.
satisfying the duality condition (\ref{dual}). Then three types of
elementary quadrilaterals of ${\Bbb Z}^3$ have the following
cross-ratios:
\begin{eqnarray}            \label{cross-ratio_delta/delta}
q(z_{k,\ell,m},z_{k+1,\ell,m},z_{k+1,\ell,m-1},z_{k,\ell,m-1})
&=&\dfrac{\Delta_1}{\Delta_3},\nonumber\\
q(z_{k,\ell,m},z_{k,\ell,m+1},z_{k,\ell-1,m+1},z_{k,\ell-1,m})
&=&\dfrac{\Delta_3}{\Delta_2},\\
q(z_{k,\ell,m},z_{k,\ell+1,m},z_{k-1,\ell+1,m},z_{k-1,\ell,m})
&=&\dfrac{\Delta_2}{\Delta_1},\nonumber
\end{eqnarray}
Same identities hold with $z$ replaced by $z^*$.

(ii) Given a solution $z:{\Bbb Z}^3\to{\Bbb C}$ to
(\ref{cross-ratio_delta/delta}) formulas (\ref{dual}) determine
(uniquely up to translation) a solution $z^*:{\Bbb Z}^3\to{\Bbb
C}$ of the same system (\ref{cross-ratio_delta/delta}).
\end{theorem}

It is obvious that the zero curvature condition (\ref{zero curv
cond}) is fulfilled for every closed contour in ${\Bbb Z}^3$ if
and only if it holds for all elementary quadrilaterals. It is easy
to check that the transition matrices $L^{(n)}$ defined above
satisfy the zero curvature condition.
\begin{theorem}                     \label{t.lax}
Let $z,z^*:{\Bbb Z}^3\to{\Bbb C}$ be a solution to (\ref{dual})
and $\ee_1, \ee_2, \ee_3, \ee_4$ be consecutive positively
oriented edges of an elementary quadrilateral of $\cQ\cL$. Then
the zero curvature condition
$$
L(\ee_4, \lambda)L(\ee_3, \lambda)L(\ee_2, \lambda)L(\ee_1,
\lambda)=I
$$
holds with $L(\ee,\lambda)$ defined by (\ref{L}). Moreover let the
elements
\begin{equation*}
L^{(n)}(\ee,\lambda)=(1-\lambda^2
\Delta_n)^{-1/2}\left(\begin{array}{cc}
1 & \lambda f\\
\lambda g & 1 \end{array}\right), \qquad fg=\Delta_n, n=1,2,3
\end{equation*}
of $G[\lambda]$ be attached to oriented edges
$\ee=(p_{out},p_{in})\in E_n$. Then the zero curvature condition
on ${\Bbb Z}^3$ is equivalent to existence of $z,z^*:{\Bbb
Z}^3\to{\Bbb C}$ such that the factorization
$$
f(\ee)=z(p_{in})-z(p_{out}),\qquad g(\ee)=z^*(p_{in})-z^*(p_{out})
$$
holds. So defined $z,z^*$ satisfy (\ref{dual},
\ref{cross-ratio_delta/delta}).
\end{theorem}

The zero curvature condition in Theorem \ref{t.lax} is a
generalization of the Lax pair found in \cite{NC} for the discrete
conformal mappings \cite{BP1}. The zero curvature condition
implies the existence of the wave function $\Psi:{\Bbb Z}^3\to
G[\lambda]$. The last one can be used to restore the fields $z$
and $z^*$. There holds the following result having many analogs in
the differential geometry described by integrable systems ("Sym
formula", see, e.q., \cite{BP2}).
\begin{theorem}                         \label{t.sym}
Let $\Psi(p, \lambda)$ be the solution of (\ref{wave function
evol}) with the initial condition $\Psi(p=0,\lambda)=I$. Then the
fields $z,z^*$ may be found as
\begin{equation}\label{Sym}
\left.\frac{d\Psi_{k,\ell,m}}{d\lambda}\right|_{\lambda=0}=
\left(\begin{array}{cc} 0 & z_{k,\ell,m}-z_{0,0,0}\\
z^*_{k,\ell,m}-z^*_{0,0,0}& 0
\end{array} \right).
\end{equation}
\end{theorem}
This simple observation turns out to be useful for analytic
constructions of solutions in particular in Section \ref{s.z^c}.

Interpreting the lattice $\cQ\cL$ as $\{(k,\ell,m)\in{\Bbb Z}^3 :
\mid k+\ell +m  \mid\leq 1\}$ one arrives at the following
\begin{corollary}
Theorems \ref{t.cross-ratio delta/delta}, \ref{t.lax}, \ref{t.sym}
hold for the lattice $\cQ\cL$, i.e. if in the statements one
replaces ${\Bbb Z}^3$ by $\cQ\cL$, $p$ by $\zz$, $E_n$ by $E^Q_n$
and assumes $k+\ell +m=0$ in (\ref{cross-ratio_delta/delta}).
\end{corollary}

Let us return to hexagonal circle patterns and explain their
relation to the discrete integrable systems of this Section. To
obtain the lattice $\cQ\cL$ we extend the intersection points of a
hexagonal circle pattern by the conformal centers of the circles.
The image of so defined mapping $z:V(\cQ\cL)\to \hat{\Bbb C}$
consists of the intersection points $z(V_1\cup V_{-1})$ and the
conformal centers $z(V_0)$. The edges connecting the points on the
circles with their centers correspond to the edges of the
quadrilateral lattice $\cQ\cL$. Whereas the angles $\alpha_n$ are
associated to three types $E^H_n,n=1,2,3$ of the edges of the
hexagonal lattice, constants $\delta_n$ defined in
(\ref{eq.delta}) are associated to three types $E^Q_n,n=1,2,3$ of
the edges of the quadrilateral lattice. Identifying them with
\begin{equation}                        \label{eq.Deltadelta}
\Delta_n=e^{-i\delta_n}
\end{equation}
in the matrices $L^{(n)}(\lambda)$ above one obtains a zero
curvature representation for hexagonal circle patterns with
constant angles.

\begin{theorem}
The intersection points and the conformal centers of the circles
$z:\cQ\cL\to \hat{\Bbb C}$ of a hexagonal circle pattern with
constant angles $\alpha_n, n=1,2,3$ satisfy the cross-ratio system
(\ref{cross-ratio_delta/delta}) on the lattice $\cQ\cL$ with
$\Delta_n, n=1,2,3$ determined by (\ref{eq.delta},
\ref{eq.Deltadelta}).
\end{theorem}
The theorem follows from the identification of
(\ref{cross-ratio_delta/delta}) with (\ref{cross-ratio=angle})
using (\ref{eq.delta}) and (\ref{eq.Deltadelta}).

The duality transformation introduced above for arbitrary mappings
satisfying the cross-ratio equations preserves the class of such
mappings coming from circle patterns with constant angles if one
extends the intersection points by the Euclidean centers of the
circles. The last one is nothing but the dual circle pattern of
Definition \ref{d.dual}.

\begin{theorem}
Let $z:\cQ\cL\to \hat{\Bbb C}$ be a hexagonal circle pattern $CP$
with constant angles extended by the Euclidean centers of the
circles. Then the dual circle pattern $CP^*$ together with the
Euclidean centers of the circles is given by the dual mapping
$z^*:\cQ\cL\to \hat{\Bbb C}$.
\end{theorem}
Since $\Delta_n, n=1,2,3$ are unitary (\ref{eq.Deltadelta}),
relation (\ref{dual.radii}) follows directly from (\ref{dual}).
The patterns have the same intersection angles since the
cross-ratio equation (\ref{cross-ratio_delta/delta}) for $z$ and
$z^*$ coincide.

\section{$z^c$ and $\log z$ patterns}           \label{s.z^c}

To construct hexagonal circle patterns analogs of holomorphic
functions $z^c$ and $\log z$ we use the analytic point description
presented in Section \ref{s.DualLax}. Recall that we interpret the
lattice $\cQ\cL$ of the intersection points and the centers of the
circles as a subset $\mid k+\ell +m\mid\leq 1$ of ${\Bbb Z}^3$ and
thus study the zero-curvature representation and $\Psi$-function
on ${\Bbb Z}^3$.

A fundamental role in the presentation of this Section is played
by a {\em non-autonomous constraint} for the solutions of the
cross-ratio system (\ref{cross-ratio_delta/delta})
\begin{eqnarray}                            \label{constraint}
bz_{k,\ell,m}^2+cz_{k,\ell,m}+d=
&2(k-a_1)\dfrac{(z_{k+1,\ell,m}-z_{k,\ell,m})
(z_{k,\ell,m}-z_{k-1,\ell,m})}{z_{k+1,\ell,m}-z_{k-1,\ell,m}}+\nonumber\\
&2(\ell-a_2)\dfrac{(z_{k,\ell+1,m}-z_{k,\ell,m})
(z_{k,\ell,m}-z_{k,\ell-1,m})}{z_{k,\ell+1,m}-z_{k,\ell-1,m}}+\\
&2(m-a_3)\dfrac{(z_{k,\ell,m+1}-z_{k,\ell,m})
(z_{k,\ell,m}-z_{k,\ell,m-1})}{z_{k,\ell,m+1}-z_{k,\ell,m-1}}.\nonumber
\end{eqnarray}
where $b,c,d,a_1,a_2,a_3\in {\Bbb C}$ are arbitrary. Note that the
form of the constraint is invariant with respect to M\"obius
transformations.

Our presentation in this section consists of three parts. First,
we explain the origin of the constraint (\ref{constraint})
deriving it in the context of isomonodromic solutions of
integrable systems. Then we show that it is compatible with the
cross-ratio system (\ref{cross-ratio_delta/delta}), i.e. there
exist non-trivial solutions of
(\ref{cross-ratio_delta/delta},\ref{constraint}). And finally we
specify parameters of these solutions to obtain circle pattern
analogs of holomorphic mappings $z^c$ and $\log z$.

We chose a different gauge of the transition matrices to simplify
formulas. Let us orient the edges of ${\Bbb Z}^3$ in the
directions of increasing $k+\ell+m$. We conjugate
$L^{(n)}(\lambda)$ of positively oriented edges with the matrix
${\rm diag}(1,\lambda)$, and then multiply by
$(1-\Delta_n\lambda^2)^{1/2}$ in order to get rid of the
normalization of the determinant. Writing then $\mu$ for
$\lambda^2$, we end up with the matrices
\begin{equation}
\cL^{(n)}(\ee,\mu)=
\left(\begin{array}{cc} 1 & z_{in}-z_{out}\\
\mu \dfrac{\Delta_n}{z_{in}-z_{out}} & 1\end{array}\right).
\end{equation}
associated to the edge $\ee=(p_{out},p_{in})\in E_n$ oriented in
the direction of increasing of $k+\ell+m$. Each elementary
quadrilateral of ${\Bbb Z}^3$ has two consecutive positively
oriented pairs of edges $\ee_1, \ee_2$ and $\ee_3, \ee_4$.  The
zero-curvature condition turns into
$$
\cL^{(n_1)}(\ee_2)\cL^{(n_2)}(\ee_1)=
\cL^{(n_2)}(\ee_4)\cL^{(n_1)}(\ee_3).
$$
Then the values of the wave function $\Phi$ in neighboring
vertices are related by the formulas
\begin{equation}                    \label{wave evolution in mu}
\left\{\begin{array}{l}
\Phi_{k+1,\ell,m}(\mu)=\cL^{(1)}(\ee,\mu)\Phi_{k,\ell,m}(\mu),\quad
\ee\in E_1,\\
\Phi_{k,\ell+1,m}(\mu)=\cL^{(2)}(\ee,\mu)\Phi_{k,\ell,m}(\mu),\quad
\ee\in E_2,\\
\Phi_{k,\ell,m+1}(\mu)=\cL^{(3)}(\ee,\mu)\Phi_{k,\ell,m}(\mu),
\quad \ee\in E_3.
\end{array}\right.
\end{equation}
We call a solution $z:{\Bbb Z}^3\mapsto{\Bbb C}$ of the equations
(\ref{cross-ratio_delta/delta}) {\em isomonodromic} (cf.
\cite{I}), if there exists a wave function $\Phi:{\Bbb Z}^3\mapsto
{\rm GL}(2,{\Bbb C})[\mu]$ satisfying (\ref{wave evolution in mu})
and some linear differential equation in $\mu$:
\begin{equation}                                \label{eq in mu}
\frac{d}{d\mu}\Phi_{k,\ell,m}(\mu)=\cA_{k,\ell,m}(\mu)\Phi_{k,\ell,m}(\mu),
\end{equation}
where $\cA_{k,\ell,m}(\mu)$ are $2\times 2$ matrices, meromorphic
in $\mu$, with the poles whose position and order do not depend on
$k,\ell,m$.

It turns out that the simplest non-trivial isomonodromic solutions
satisfy the constraint (\ref{constraint}). Indeed, since $\det
\cL^{(n)}(\mu)$ vanishes at $\mu=1/\Delta_n$ the logarithmic
derivative of $\Phi(\mu)$ must be singular in these points. We
assume that these singularities are as simple as possible, i.e.
simple poles.

\begin{theorem}                                 \label{mu-theorem}
Let $z:{\Bbb Z}^3\to {\Bbb C}$ be an isomonodromic solution to
(\ref{cross-ratio_delta/delta}) with the matrix $\cA_{k,\ell,m}$
in (\ref{eq in mu}) of the form
\begin{equation}                            \label{A-matrixStructure}
\cA_{k,\ell,m}(\mu)=\dfrac{C_{k,\ell,m}}{\mu}+
\sum_{n=1}^3\dfrac{B^{(n)}_{k,\ell,m}}{\mu -\dfrac{1}{\Delta_n}}
\end{equation}
with $\mu$-independent matrices $C_{k,\ell,m},B^{(n)}_{k,\ell,m}$
and normalized\footnote{As explained in the proof of the theorem
this normalization can be achieved without loss of generality.}
trace ${\rm tr}\cA_{0,0,0}(\mu)=0$. Then these matrices have the
following form:
\begin{eqnarray*}
C_{k,\ell,m} &=&\dfrac{1}{2} \left(\begin{array}{cc}
-bz_{k,\ell,m}-c/2
& bz_{k,\ell,m}^2+cz_{k,\ell,m}+d\\
b & bz_{k,\ell,m}+c/2\end{array}\right)\\
B^{(1)}_{k,\ell,m} &=&\dfrac{k-a_1}{z_{k+1,\ell,m}-z_{k-1,\ell,m}}
\left(\begin{array}{cc} z_{k+1,\ell,m}-z_{k,\ell,m}
& (z_{k+1,\ell,m}-z_{k,\ell,m})(z_{k,\ell,m}-z_{k-1,\ell,m})\\
1 & z_{k,\ell,m}-z_{k-1,\ell,m}\end{array}\right)+
\dfrac{a_1}{2}I\\
B^{(2)}_{k,\ell,m}
&=&\dfrac{\ell-a_2}{z_{k,\ell+1,m}-z_{k,\ell-1,m}}
\left(\begin{array}{cc} z_{k,\ell+1,m}-z_{k,\ell,m}
& (z_{k,\ell+1,m}-z_{k,\ell,m})(z_{k,\ell,m}-z_{k,\ell-1,m})\\
1 & z_{k,\ell,m}-z_{k,\ell-1,m}\end{array}\right)+
\dfrac{a_2}{2}I\\
B^{(3)}_{k,\ell,m} &=&\dfrac{m-a_3}{z_{k,\ell,m+1}-z_{k,\ell,m-1}}
\left(\begin{array}{cc} z_{k,\ell,m+1}-z_{k,\ell,m}
& (z_{k,\ell,m+1}-z_{k,\ell,m})(z_{k,\ell,m}-z_{k,\ell,m-1})\\
1 & z_{k,\ell,m}-z_{k,\ell,m-1}\end{array}\right)+ \dfrac{a_3}{2}I
\end{eqnarray*}
and $z_{k,\ell,m}$ satisfies (\ref{constraint}).

Conversely, any solution $z:{\Bbb Z}^3\to {\Bbb C}$ to the system
(\ref{cross-ratio_delta/delta},\ref{constraint}) is isomonodromic
with $\cA_{k,\ell,m}(\mu)$ given by the formulas above.
\end{theorem}
The proofs of this theorem as well as of the next one are
presented in Appendix \ref{Appendix.lambda}. We prove Theorem
\ref{mu-theorem} by computing the compatibility conditions of
(\ref{wave evolution in mu}) and (\ref{eq in mu}) and Theorem
\ref{t.compatibility} by showing the solvability of a reasonably
posed Cauchy problem.

\begin{theorem}                                 \label{t.compatibility}
For arbitrary $b,c,d,a_1,a_2,a_3\in {\Bbb C}$ the constraint
(\ref{constraint}) is compatible with the cross-ratio equations
(\ref{cross-ratio_delta/delta}), i.e. there are non-trivial
solutions of (\ref{cross-ratio_delta/delta},\ref{constraint}).
\end{theorem}


Further we deal with the special case of (\ref{constraint}) where
$$
b=a_1=a_2=a_3=0.
$$
If $c\not=0$ one can always assume $d=0$ in (\ref{constraint}) by
shifting $z\to z-\dfrac{d}{c}$:
\begin{eqnarray}                            \label{constraint.z^c}
cz_{k,\ell,m}= &2k\dfrac{(z_{k+1,\ell,m}-z_{k,\ell,m})
(z_{k,\ell,m}-z_{k-1,\ell,m})}{z_{k+1,\ell,m}-z_{k-1,\ell,m}}+\nonumber\\
&2\ell\dfrac{(z_{k,\ell+1,m}-z_{k,\ell,m})
(z_{k,\ell,m}-z_{k,\ell-1,m})}{z_{k,\ell+1,m}-z_{k,\ell-1,m}}+\\
&2m\dfrac{(z_{k,\ell,m+1}-z_{k,\ell,m})
(z_{k,\ell,m}-z_{k,\ell,m-1})}{z_{k,\ell,m+1}-z_{k,\ell,m-1}}.\nonumber
\end{eqnarray}
To define a circle pattern analogue of $z^c$ it is natural to
restrict the mapping to the following two subsets of ${\Bbb Z}^3$:
$$
Q=\{(k,\ell,m)\in {\Bbb Z}^3\mid k\geq 0, \ell\geq 0, m\leq 0\},\
H=\{(k,\ell,m)\in {\Bbb Z}^3\mid m\leq 0\}\subset {\Bbb Z}^3.
$$
\begin{corollary}
The solution $z:Q\to{\Bbb C}$ of the system
(\ref{cross-ratio_delta/delta}) satisfying the constraint
(\ref{constraint.z^c}) is uniquely determined by its values
\begin{equation}                    \label{cauchy.z}
z_{1,0,0},\ z_{0,1,0},\ z_{0,0,-1}.
\end{equation}
\end{corollary}
{\em Proof.} Using the constraint one determines the values along
the coordinate lines $z_{n,0,0}, z_{0,n,0}, z_{0,0,-n}$  $\forall
n\in{\Bbb N}$. Then all other $z_{k,\ell,m},\ k,\ell,-m\in{\Bbb
N}$ in consecutive order are determined through the cross-ratios
(\ref{cross-ratio_delta/delta}). Computations using different
cross-ratios give the same result due to the following lemma about
the eighth point.
\begin{lemma}                                   \label{hexagon.lemma}
Let $H(p_1+1/2,p_2+1/2,p_3-1/2), p=(p_1,p_2,p_3)\in{\Bbb Z}^3$ be
the elementary hexahedron (lying in ${\Bbb C}$) with the vertices
$z_{p+(k,\ell,-m)},\ k,\ell,m\in\{0,1\}$ and let the cross-ratios
of the opposite cites of $H(p_1+1/2,p_2+1/2,p_3-1/2)$ be equal
\begin{eqnarray*}
q(z_p, z_{p+(0,0,-1)}, z_{p+(0,1,-1)}, z_{p+(0,1,0)} )=
q(z_{p+(1,0,0)},z_{p+(1,0,-1)}, z_{p+(1,1,-1)},  z_{p+(1,1,0)})&=:&q_1\\
q(z_p, z_{p+(1,0,0)}, z_{p+(1,0,-1)},z_{p+(0,0,-1)} )=
q(z_{p+(0,1,0)},z_{p+(1,1,0)} , z_{p+(1,1,-1)},z_{p+(0,1,-1)} )&=:&q_2\\
q(z_p, z_{p+(0,1,0)}, z_{p+(1,1,0)}, z_{p+(1,0,0)} )=
q(z_{p+(0,0,-1)},z_{p+(0,1,-1)}, z_{p+(1,1,-1)},z_{p+(1,0,-1)}
)&=:&q_3
\end{eqnarray*}
with $q_1 q_2 q_3=1$. Then all the vertices of the hexahedron are
uniquely determined through four given $z_p$, $z_{p+(1,0,0)}$,
$z_{p+(0,1,0)}$, $z_{p+(0,0,-1)}$.
\end{lemma}
The relation of solutions of
(\ref{cross-ratio_delta/delta},\ref{constraint.z^c}) to circle
patterns is established in the following
\begin{theorem}                   \label{t.Z^3circlepattern}
The solution $z:Q\to{\Bbb C}$ of the system
(\ref{cross-ratio_delta/delta},\ref{constraint.z^c}) with the
initial data
\begin{equation}                    \label{cauchy.betagamma}
z_{1,0,0}=1,\ z_{0,1,0}=e^{i\beta},\ z_{0,0,-1}=e^{i\gamma}
\end{equation}
and unitary cross-ratios $q_n=e^{-2i\alpha_n}$ determines a circle
pattern. For all $(k,\ell,m)\in Q$ with even $k+\ell+m $ the
points $z_{k\pm 1,\ell,m}, z_{k,\ell\pm 1,m}, z_{k,\ell,m\pm 1}$
lie on a circle with the center $z_{k,\ell,m}$.
\end{theorem}
{\it Proof}. We say that a quadrilateral is of the {\em kite form}
if it has two pairs of equal adjacent edges, (and that a
hexahedron is of the kite form if all its cites are of the kite
form). These quadrilaterals have unitary cross-ratios with the
argument equal to the angle between the edges (see Lemma
\ref{l.cross-ratio=angle}). Our proof of the theorem is based on
two simple observations:
\begin{itemize}
  \item[(i)] If a quadrilateral has a pair of adjacent edges of equal length  and
  unitary cross-ratio then it is of the kite form,
  \item[(ii)]If the cross-ratio of a quadrilateral is equal to
  $q(z_1,z_2,z_3,z_4)=e^{-2i\alpha}$ where $\alpha$ is the angle between
  the edges $(z_1,z_2)$ and $(z_2,z_3)$
  (as in Figure \ref{fig:cross-ratio}) then the quadrilateral
  is of the kite form.
\end{itemize}
Applying these observation to the elementary hexahedron in Lemma
\ref{hexagon.lemma} one obtains
\begin{lemma}                                   \label{hexagon-kite.lemma}
Let three adjacent edges of the elementary hexahedron in Lemma
\ref{hexagon.lemma} have equal length and the cross-ratios of the
cites be unitary $q_n=e^{-2i\alpha_n}$. Then the hexahedron is of
the kite form.
\end{lemma}

The constraint (\ref{constraint.z^c}) with $\ell=m=0$ implies by
induction
$$
\mid z_{2n + 1,0,0}-z_{2n,0,0}\mid = \mid
z_{2n,0,0}-z_{2n-1,0,0}\mid,
$$
and all the points $z_{n,0,0}$ lie on the real axis. The same
equidistance and straight line properties hold true for $\ell$ and
$m$ axes. Also by induction one shows that all elementary
hexahedra of the mapping $z:Q\to{\Bbb C}$ are of the kite form.
Indeed, due to Lemma \ref{hexagon-kite.lemma} the initial
conditions (\ref{cauchy.betagamma}) imply that the hexahedron
$H(1/2,1/2,-1/2)$ is of the kite form. Since the axis vertices lie
on the straight lines,  we get, for example, that the angle
between the edges $(z_{1,1,0},z_{1,0,0})$ and
$(z_{1,0,0},z_{2,0,0})$ is $\alpha_3$. Now applying the
observation (ii) to the cites of $H(3/2,1/2,-1/2)$ we see that
$$
\mid z_{2,0,0}-z_{1,0,0}\mid=\mid z_{2,0,0}-z_{2,1,0}\mid=\mid
z_{2,0,0}-z_{2,0,-1}\mid.
$$
Lemma \ref{hexagon-kite.lemma} implies $H(3/2,1/2,-1/2)$ is of the
kite form. Proceeding further this way and controlling the lengths
of the edges meeting at $z_{k,\ell,m}$ with even $k+\ell+m$ and
the angles at $z_{k,\ell,m}$ with odd $k+\ell+m$ one proves the
statement for all hexahedra. \vspace{3mm}

Denote the vertices of the hexagonal grid by (see Figure
\ref{fig:QLattice})
\begin{eqnarray*}
Q_H&=&\{(k,\ell,m)\in Q:\ \mid k+\ell+m\mid \leq 1\},\\
H_H&=&\{(k,\ell,m)\in H:\ \mid k+\ell+m\mid \leq 1\}.
\end{eqnarray*}
The half-plane $H_H$ consists of the sector $Q_H$ and its two
images under the rotations with the angles $\pm \pi/3$.
\begin{corollary}
The mapping $z:H_H \to {\Bbb C}$ given by
(\ref{cross-ratio_delta/delta},\ref{constraint.z^c}) with unitary
cross-ratios $q_n=e^{-2i\alpha_n}$ and unitary initial data
$$
\mid z_{\pm 1,0,0}\mid=\mid z_{0,\pm 1,0}\mid= \mid
z_{0,0,0}\mid=1
$$
is a hexagonal circle pattern with constant angles $\alpha_n$
extended by the Euclidean centers of the circles. The centers of
the corresponding circles are the images of the points with
$k+\ell+m=0$.
\end{corollary}

\begin{figure}[htbp]
  \begin{center}
\centerline{\includegraphics[width=0.45\hsize]{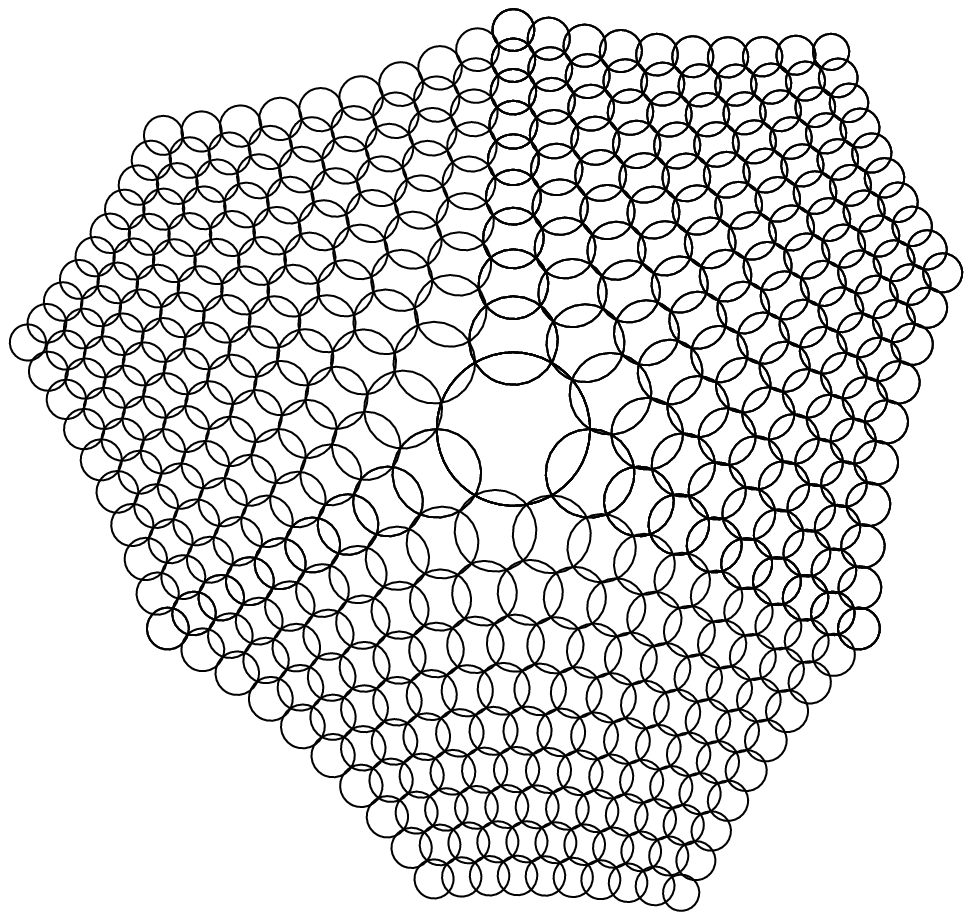}%
\includegraphics[width=0.45\hsize]{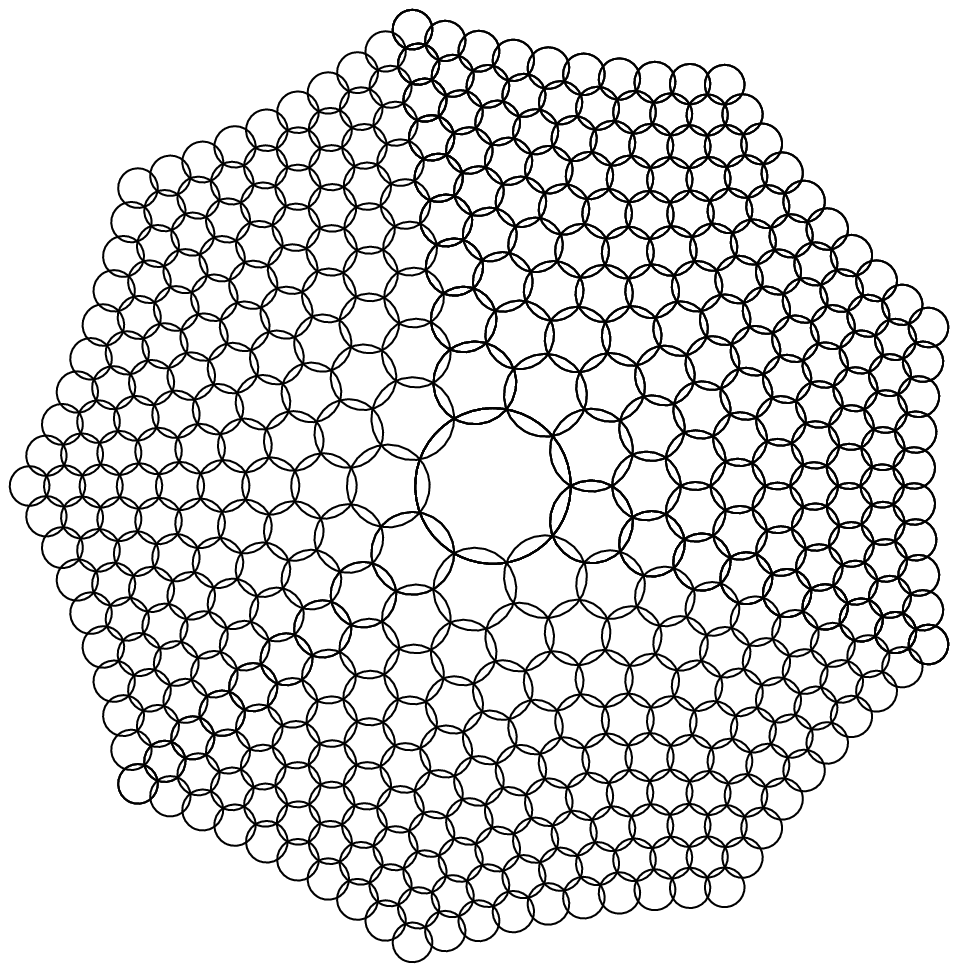}}\caption{Non-isotropic
and isotropic circle patterns $z^{2/3}$.} \label{fig:z^c.general}
  \end{center}
\end{figure}

The circle patterns determined by most of the initial data
$\beta,\gamma\in{\Bbb R}$ are quite irregular. But for a special
choice of these parameters one obtains a regular circle pattern
which we call the hexagonal circle pattern $z^c$ motivated by the
asymptotic of the constraint (\ref{constraint.z^c}) as
$k,\ell,m\to\infty$.
\begin{definition}                      \label{def.z^c}
The hexagonal circle pattern $z^c, \ 0<c<2$ (extended by the
Euclidean centers of the circles) is the solution $z:H_H \to{\Bbb
C}$ of (\ref{cross-ratio_delta/delta},\ref{constraint.z^c}) with
$$
\dfrac{\Delta_{n+1}}{\Delta_{n+2}}=e^{2i\alpha_n}, \qquad n \pmod
3,
$$
and the initial conditions
\begin{equation}                    \label{cauchy.z^c.general}
z_{1,0,0}=1,\ z_{0,1,0}=e^{ic(\alpha_1+\alpha_2)},\
z_{0,0,-1}=e^{ic\alpha_2},\ z_{-1,0,0}=e^{ic\pi},\
z_{0,-1,0}=e^{-ic\alpha_3}.
\end{equation}
\end{definition}
Motivated by our computer experiments (see in particular Figure
\ref{fig:z^c.general}) and the corresponding result for Schramm's
circle patterns with the combinatorics of the square grid
\cite{AB,A} we conjecture that the hexagonal $z^c$ is embedded,
i.e. the interiors of all elementary quadrilaterals
$(z(\zz_1),z(\zz_2),z(\zz_3),z(\zz_4)), \mid \zz_{i+1}-\zz_i\mid
=1, i \pmod 4$ are disjoint.

In the isotropic case when all the intersection angles are the
same
$$
\alpha_1=\alpha_2=\alpha_3=\dfrac{\pi}{3}.
$$
one has
\begin{equation}                        \label{Delta_123}
\Delta_1=\omega \Delta_2=\omega^2\Delta_3,
\end{equation}
and due to the symmetry it is enough to restrict the mapping to
$Q_H$. The initial conditions (\ref{cauchy.z^c.general}) become
\begin{equation}                    \label{cauchy.z^c}
z_{1,0,0}=1,\ z_{0,1,0}=e^{2\pi ic/3},\ z_{0,0,-1}=e^{\pi ic/3}.
\end{equation}

Like for the smooth $z^c$ the images of the coordinate axes $\arg
\zz=\dfrac{\pi}{6}, \dfrac{\pi}{3}, \dfrac{\pi}{2}$ are the axes
with the arguments $\dfrac{\pi c}{6}, \dfrac{\pi c}{3}$ and
$\dfrac{\pi c}{2}$ respectively:
$$
\arg z_{n,0,-n}=\dfrac{\pi c}{6},\ \arg z_{0,n,-n}=\dfrac{\pi
c}{2},\ \arg z_{n,n,m}=\dfrac{\pi c}{3},\ m=-2n; -2n\pm 1, \qquad
\forall n\in{\Bbb N}.
$$

\begin{figure}[htbp]
  \begin{center}
\includegraphics[width=0.3\hsize]{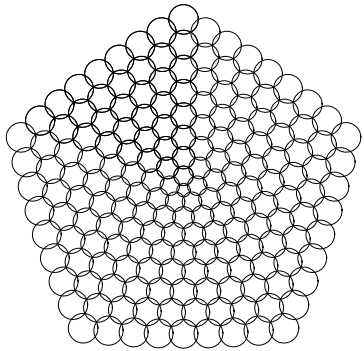}
\caption{Circle pattern $z^{6/5}$.} \label{fig:z65}
  \end{center}
\end{figure}

For $c=\dfrac{6}{q},\ q\in{\Bbb N}$ the circle pattern $z^c(Q_H)$
together with its additional $q-1$ copies rotated by the angles
$2\pi n/q,\ n=1,\ldots,q-1$ comprise a circle pattern covering the
whole complex plane. This pattern is hexagonal at all circles
beside the central one which has $q$ neighboring circles. The
circle pattern $z^{6/5}$ is shown in Figure \ref{fig:z65}.

Definition \ref{def.z^c} is given for $0<c<2$. For $c=2$ the radii
of the circles intersecting the central circle of the pattern
become infinite. To get finite radii the central circle should
degenerate to a point
\begin{equation}                            \label{z^2.a}
z_{0,0,0}=z_{1,0,0}=z_{0,0,-1}=z_{0,1,0}=0.
\end{equation}
The mapping $z: Q\to {\Bbb C}$ is then uniquely determined by
$z_{2,0,0}, z_{0,2,0}, z_{0,0,-2}, z_{1,1,0}, z_{1,0,-1},
z_{0,1,-1}$, the values of which can be derived taking the limit
$c\to 2-0$. Indeed, consider the quader $m=0$ with the
cross-ratios of all elementary quadrilaterals equal to
$q=e^{-2i\alpha}$. Choosing $z_{1,0,0}=\epsilon,
z_{0,1,0}=\epsilon e^{ic\alpha}$ as above we get
$z_{2,0,0}=\dfrac{2\epsilon}{2-c}$. Normalization
$$
c=2-2\epsilon
$$
yields
\begin{equation}                            \label{z_200}
z_{2,0,0}=1,\quad z_{0,2,0}=e^{ic\alpha}.
\end{equation}
Since the angles of the triangle with the vertices
$z_{0,0,0},z_{1,0,0},z_{1,1,0}$ are $\dfrac{c\alpha}{2},
\pi-\alpha, \alpha\epsilon$ respectively, one obtains
$$
z_{1,1,0}=\epsilon+R_\epsilon e^{i\alpha},\qquad
R_\epsilon=\epsilon\dfrac{\sin (c\alpha/2)}{\sin
(\alpha\epsilon)}.
$$
In the limit $\epsilon\to 0$ we have
\begin{equation}                            \label{z_110}
z_{1,1,0}=\dfrac{\sin\alpha}{\alpha}e^{i\alpha}.
\end{equation}
Observe that in our previous notations
$\alpha=\alpha_1+\alpha_2=\pi-\alpha_3$. Finally, the same
computations for the sectors $\{\ell=0, k\geq 0, m\leq 0\}$ and
$\{k=0, \ell\geq 0, m\leq 0\}$ provide us with the following data:
\begin{eqnarray}                            \label{z^2.b}
&z_{2,0,0}=1,\ z_{0,0,-2}=e^{2i\alpha_2},\ z_{0,2,0}=e^{2i(\alpha_1+\alpha_2)},\nonumber\\
&z_{1,0,-1}=\dfrac{\sin\alpha_2}{\alpha_2}e^{i\alpha_2},\
z_{0,1,-1}=\dfrac{\sin\alpha_1}{\alpha_1}e^{i(\alpha_1+2\alpha_2)},\
z_{1,1,0}=\dfrac{\sin(\alpha_1+\alpha_2)}{\alpha_1+\alpha_2}e^{i(\alpha_1+\alpha_2)}.
\end{eqnarray}
In the same way the initial data for other two sectors of $H_H$
are specified. The hexagonal circle pattern with these initial
data and $c=2$ is an analog of the holomorphic mapping $z^2$ .

The duality transformations preserves the class of circle patterns
we defined.
\begin{theorem}                                 \label{c^*}
$$
(z^c)^*=z^{c^*},\qquad c^*={2-c},
$$
where $(z^c)^*$ is the hexagonal circle pattern dual to the circle
pattern $z^c$ and normalized to vanish at the origin
$(z^c)^*(\zz=0)=0$.
\end{theorem}
{\it Proof.} Let us consider the mapping on the whole $Q$. It is
easy to see that on the axes the duality transformation
(\ref{dual}) preserves the form of the constraint with $c$ being
replaced by $c^*=2-c$. Then the constraint with $c^*$ holds for
all points of $Q$ due to the compatibility in Theorem
\ref{t.compatibility}. Restriction to $\mid k+\ell+m\mid\leq 1$
implies the claim. \vspace{3mm}

The smooth limit of the duality transformation of the hexagonal
patterns is the following transformation $f\mapsto f^*$ of
holomorphic functions:
$$
(f^*(z))'=\dfrac{1}{f'(z)}.
$$
The dual of $f(z)=z^2$ is, up to constant, $f^*(z)=\log z$.
Motivated by this observation, we define the hexagonal circle
pattern $\log z$ as the dual to the circle pattern $z^2$:
$$
\log z:=(z^2)^*.
$$
The corresponding constraint (\ref{constraint})
\begin{eqnarray}                            \label{constraint.log}
2k\dfrac{(z_{k+1,\ell,m}-z_{k,\ell,m})
(z_{k,\ell,m}-z_{k-1,\ell,m})}{z_{k+1,\ell,m}-z_{k-1,\ell,m}}+
2\ell\dfrac{(z_{k,\ell+1,m}-z_{k,\ell,m})
(z_{k,\ell,m}-z_{k,\ell-1,m})}{z_{k,\ell+1,m}-z_{k,\ell-1,m}}+\nonumber\\
2m\dfrac{(z_{k,\ell,m+1}-z_{k,\ell,m})
(z_{k,\ell,m}-z_{k,\ell,m-1})}{z_{k,\ell,m+1}-z_{k,\ell,m-1}}=1.
\end{eqnarray}
can be derived as a limit $c\to +0$ (see \cite{AB} for this limit
in the square grid case). The initial data for $\log z$ are dual
to the ones of $z^2$. In our model case of the quader $m=0$
factorizing $q=\Delta_1/\Delta_2$ with $\Delta_1=1/2,
\Delta_2=e^{2i\alpha}/2$ one arrives at the following data dual to
(\ref{z^2.a},\ref{z_200},\ref{z_110}):
\begin{eqnarray}                            \label{log.initial}
&z_{0,0,0}=\infty,\ z_{1,0,0}=0,\ z_{2,0,0}=\dfrac{1}{2},\nonumber\\
&z_{0,1,0}=i\alpha,\ z_{0,2,0}=\dfrac{1}{2}+i\alpha,\
z_{1,1,0}=\dfrac{\alpha}{2\sin\alpha}e^{i\alpha}.
\end{eqnarray}
In the isotropic case the circle patterns are more symmetric and
can be described as mappings of $Q_H$.
\begin{definition}
The isotropic hexagonal circle patterns $z^2$ and $\log z$ are the
mappings $Q_H\to{\Bbb C}$ with the cross-ratios of all elementary
quadrilaterals equal\footnote{The first point in the cross-ratio
is a circle center and the quadrilaterals are positively
oriented.} to $e^{-2\pi i/3}$. The mapping $z^2$ is determined by
the constraint (\ref{constraint.z^c}) with $c=2$ and the initial
data
\begin{eqnarray*}
&z_{0,0,0}=\infty,\ z_{1,0,0}=0,\ z_{0,0,-1}=0,\
z_{0,1,0}=0,\\
&z_{2,0,0}=1,\ z_{0,0,-2}=e^{2\pi i/3},\
z_{0,2,0}=e^{4\pi i/3},\\
&z_{1,0,-1}=\dfrac{3\sqrt{3}}{2\pi}e^{\pi i/3},\
z_{0,1,-1}=-\dfrac{3\sqrt{3}}{2\pi},\
z_{1,1,0}=\dfrac{3\sqrt{3}}{4\pi}e^{2\pi i/3}.
\end{eqnarray*}
For $\log z$ the corresponding constraint is
(\ref{constraint.log}) and the initial data are
\begin{eqnarray*}
&z_{0,0,0}=\infty,\ z_{1,0,0}=0,\
z_{0,0,-1}=\dfrac{\pi}{3}i, \
z_{0,1,0}=\dfrac{2\pi}{3}i,\\
&z_{2,0,0}=\dfrac{1}{2},\
z_{0,0,-2}=\dfrac{1}{2}+\dfrac{\pi}{3}i,\
z_{0,2,0}=\dfrac{1}{2}+\dfrac{2\pi}{3}i,\\
&z_{1,0,-1}=\dfrac{\pi}{6}(\dfrac{1}{\sqrt{3}}+i),\
z_{0,1,-1}=\dfrac{\pi}{6}(\dfrac{1}{\sqrt{3}}+3i),\
z_{1,1,0}=\dfrac{\pi}{3}(-\dfrac{1}{\sqrt{3}}+i).
\end{eqnarray*}
\end{definition}

\begin{figure}[htbp]
  \begin{center}
\centerline{\includegraphics[width=0.35\hsize]{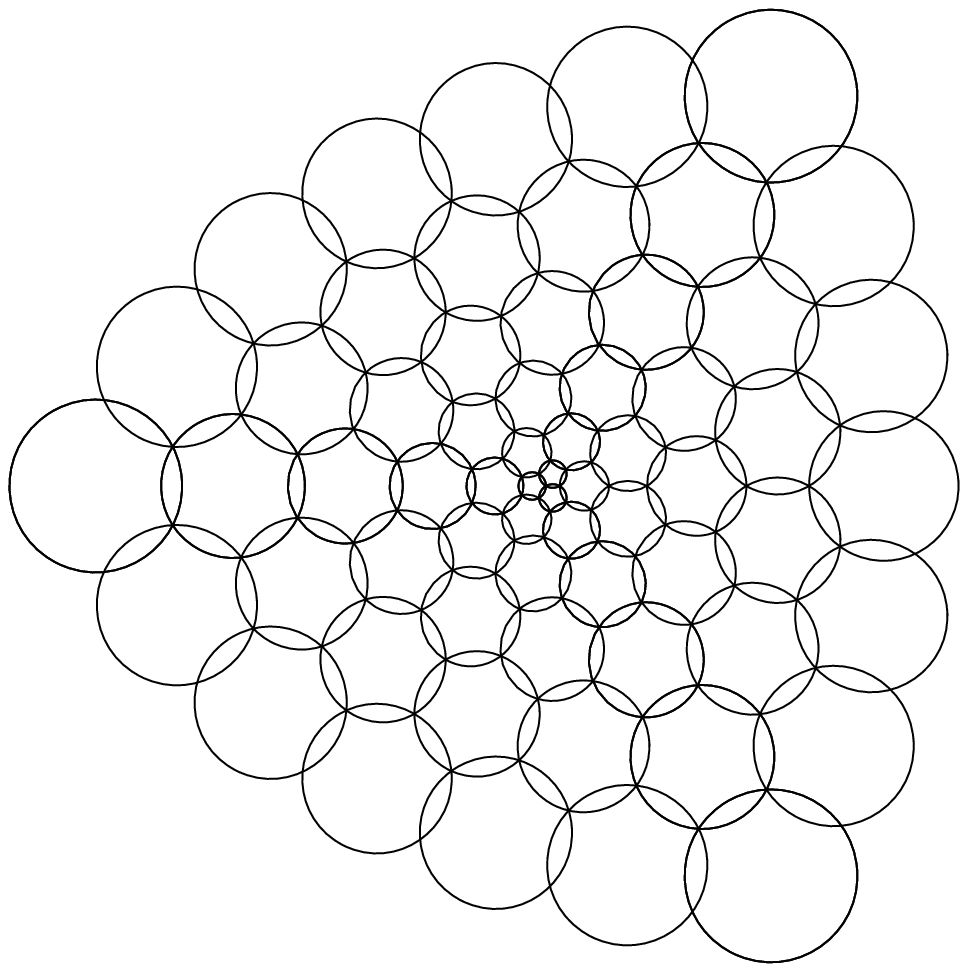}\hspace{1.0cm}%
\raise 0.55cm\hbox{\includegraphics[width=0.35\hsize]{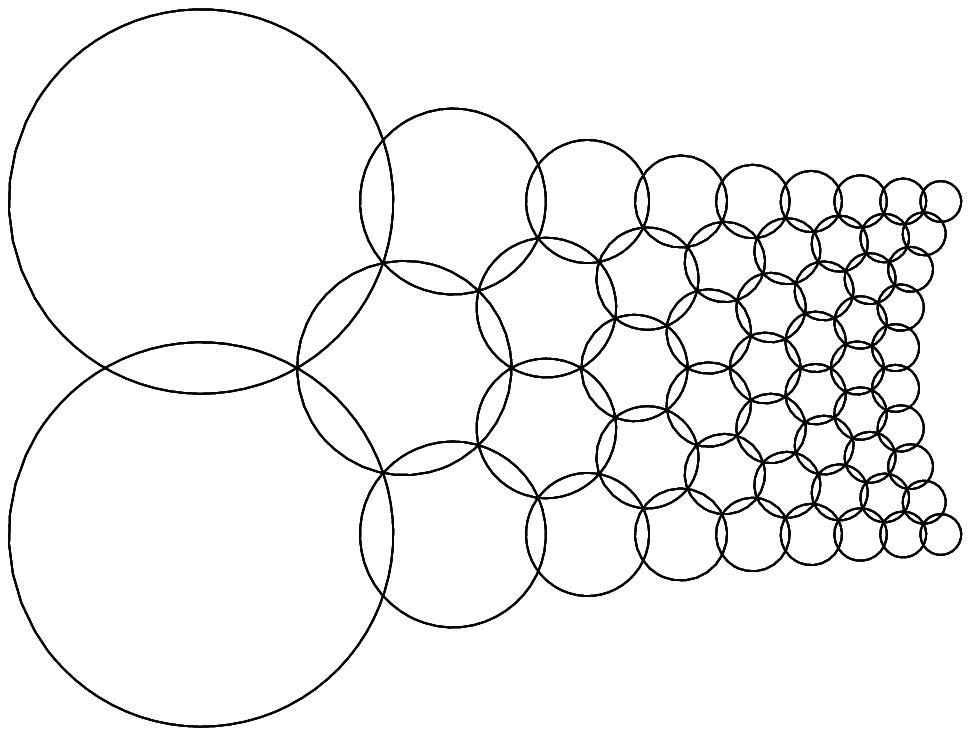}}}
\caption{Isotropic circle patterns $z^2$ and $\log
z$.} \label{fig:z^2logz}
  \end{center}
\end{figure}

The isotropic hexagonal circle patterns $z^2$ and $\log z$ are
shown in Figure \ref{fig:z^2logz}.

Starting with $z^c,\ c\in (0,2]$ one can easily define $z^c$ for
arbitrary $c$ by applying some simple transformations of hexagonal
circle patterns. The construction here is the same as for
Schramm's patterns (see Section 6 of \cite{AB} for details).
Applying the inversion of the complex plane $z\mapsto 1/z$ to the
circle pattern $z^c,\ c\in (0,2]$ one obtains a circle pattern
satisfying the constraint with $-c$. It is natural to call it the
hexagonal circle pattern $z^{-c}, \ c\in (0,2]$. Constructing the
dual circle pattern we arrive at a natural definition of
$z^{2+c}$. Intertwining the inversion and the dualization
described above, one constructs circle patterns $z^c$ for any $c$.

In particular, inverting and then dualizing
$z=k+\ell\omega+m\omega^2$ with $\Delta_1=-3, \Delta_2=-3\omega^2,
\Delta_3=-3\omega$ we obtain the circle pattern corresponding to
$z^3$:
$$
z_{k,\ell,m}= (k+\ell\omega+m\omega^2)^3-(k+\ell+m).
$$
Note that this is the central extension corresponding to
$P_\infty=0$. The points with even $k+\ell+m$ can be replaced by
the Euclidean centers of the circles. As it is shown in Section 6
of \cite{AB} the replacement $P_\infty=0$ by $P_\infty=\infty$
preserves the constraint (\ref{constraint.z^c}).

\section{Concluding remarks}
\label{s.Concluding}


We restricted ourself to the analysis of geometric and algebraic
properties of hexagonal circle patterns leaving the approximation
problem beyond the scope of this paper. The convergence of circle
patterns with the combinatorics of the square grid to the Riemann
mapping is proven by Schramm in \cite{S}. We expect that his
result can be extended to the hexagonal circle patterns defined in
this paper.

The entire circle pattern ${\rm erf} z$ also found in \cite{S}
remains to be rather mysterious. We were unable to find its
analogue in the hexagonal case, and thus, have no counter-examples
to the Doyle conjecture for hexagonal circle patterns with
constant angles. It seems that Schramm's ${\rm erf} z$ is a
feature of the square grid combinatorics.

The construction of the hexagonal circle pattern analogs of $z^c$
and $\log z$ in Section \ref{s.z^c} was based on the extension of
the corresponding integrable system to the lattice ${\Bbb Z}^3$.
Theorem \ref{t.Z^3circlepattern} claims that in this way one
obtains a circle pattern labelled by three independent indices $k,
\ell, m$. Fixing one of the indices, say $m=m_0$, one obtains a
Schramm's circle pattern with the combinatorics of the square
grid. In the same way the restriction of this "three-dimensional"
pattern to the sublattice $\mid k+\ell+m-n_0\mid\leq 1$ with some
fixed $n_0\in {\Bbb Z}$ yields a hexagonal circle pattern. In
Section \ref{s.z^c} we have defined circle patterns $z^c$ on the
sublattices $m_0=0$ and $n_0=0$. In the same way the sublattices
$m_0\neq 0$ and $n_0\neq 0$ can be interpreted as the circle
pattern analogs of the analytic function $(z+a)^c,\ a\neq 0$ with
the square and hexagonal grids combinatorics respectively.

It is unknown, whether the theory of integrable systems can be
applied to hexagonal circle packings. As we mentioned already in
the Introduction, the underdevelopment of the theory of integrable
systems on the lattices different from ${\Bbb Z}^n$ may be a
reason for this. We hope that integrable systems on the hexagonal
and related to it lattices appeared in this paper and in
\cite{BHS} as well as the the experience of integration of
hexagonal circle patterns in general will be helpful for the
progress in investigation of hexagonal circle packings.

\begin{appendix}
\section{Appendix. A Lax representation for the conformal
      description.}
\label{sec:confLax}

We will now give a Lax representation for equations
(\ref{eq:Tequation}). For $\zz\in\VHL$ let $[\zz_i,\zz]\in E^H_i$
and $m_\zz$ be the M\"obius transformation, that sends
$(\zz_1,\zz_2,\zz_3)$ to $(0,1,\infty)$. For $\ee
=[\zz,\tilde\zz]\in E^H_i$ set $L_i = m_{\tilde\zz}\circ
m^{-1}_\zz$. The $L_i$ depend on $\S^{(1)}$ and $\T_\ee$ only.
They are given in equation (\ref{eq:conformalLax}). Note that we
do not need to orient the edges since $L_i L_i^{-1}$ is the
identity if we normalize $\det L_i = 1$. The claim is that the
closing condition when multiplying the $L_i$ around one hexagon is
equivalent to equation~(\ref{eq:Tequation}):
\begin{theorem}
  Attach to each edge of $E_i$ a matrix
  $L_i(\T,\S)$ of the form:
  \begin{equation}
    \label{eq:conformalLax}
    \begin{array}{rcl}
      L_1(\T,\S) &=&
      \quadmatrix{-1}{\frac{\S-1}{\S}}%
      {\frac{\T+\S-1}{\S-1}}{1}\\[0.2cm]
      L_2(\T,\S) &=&\quadmatrix{\T\S}{1-\T\S}%
      {\S\frac{1+\T(\S-1)}{\S-1}}{-\T\S}\\[0.2cm]
      L_3(\T,\S) &=& \quadmatrix{1-\S}{-1 +\T(\frac1\S-1)
        +\S}{-\S}{\S-1}
    \end{array}
  \end{equation}
  then the zero curvature condition for each hexagon in $\FHL$
  \begin{equation}
    \label{eq:conformalZC}
    L_1(\T_4,\S) L_2(\T_5,\S) L_3(\T_6,\S) = \rho
    L_3(\T_3,\S) L_2(\T_2,\S) L_1(\T_1,\S)
  \end{equation}
  for all $\S$ is equivalent to (\ref{eq:Tequation}).
\end{theorem}
\begin{proof}
  Straight forward calculations.
\end{proof}
Equations (\ref{eq:conformalLax}) and (\ref{eq:conformalZC}) are a
Lax representation for equation (\ref{eq:Tequation}) with the
spectral parameter $\S$.

\section{Appendix. Discrete equations of Toda type and cross-ratios}
\label{Appendix.Toda}

In the following we will discuss briefly the connection between
discrete equations of Toda type and cross-ratio equations for the
square grid and the dual Kagome lattice. This interrelationship
holds in a more general setting, namely for discrete Toda systems
on graphs. This situation will be considered in a subsequent
publication.

We will start with the square grid. Let us decompose the lattice
${\Bbb Z}^2$ into two sublattices
$$
{\Bbb Z}^2_k=\{(m,n)\in {\Bbb Z}^2, (m+n\!\!\!\mod 2)=k\},\qquad
k=0,1.
$$

\begin{theorem}
Let $z: {\Bbb Z}^2\to\hat{\Bbb C}$ be a solution to the
cross-ratio equation
\begin{equation}                            \label{cross-ratio(mn)}
q(z_{m,n}, z_{m+1,n}, z_{m+1,n+1}, z_{m,n+1})=q.
\end{equation}
Then restricted to the sublattices ${\Bbb Z}^2_0$ or ${\Bbb
Z}^2_1$ it satisfies the discrete equation of Toda type
(\ref{Toda_z})
\begin{equation}                            \label{Toda(mn)}
\dfrac{1}{z_{m,n}-z_{m+1,n+1}}+\dfrac{1}{z_{m,n}-z_{m-1,n-1}}=
\dfrac{1}{z_{m,n}-z_{m+1,n-1}}+\dfrac{1}{z_{m,n}-z_{m-1,n+1}}.
\end{equation}
Conversely, given a solution $z: {\Bbb Z}^2_0\to\hat{\Bbb C}$ of
 equation (\ref{Toda(mn)}) and arbitrary $q, w\in
\hat{\Bbb C}$ there exists unique extension of $z$ to the whole
lattice $z: {\Bbb Z}^2\to\hat{\Bbb C}$ with $z_{1,0}=w$ satisfying
the cross-ratio condition (\ref{cross-ratio(mn)}). Moreover,
restricted to the other sublattice ${\Bbb Z}^2_1$ the so defined
mapping $z: {\Bbb Z}^2_1\to\hat{\Bbb C}$ satisfies the same
discrete equation of Toda type (\ref{Toda(mn)}).
\end{theorem}
\begin{proof}
  For any four complex numbers the cross-ratio equation can be written
  in a simple fraction form:
  \begin{equation}
    \label{eq:a:crAForm}
    q(u_1,u_2,u_3,u_4) = q \Leftrightarrow \frac1{u_2-u_1} -
    \frac{q}{u_2-u_3} + \frac{q-1}{u_2-u_4} = 0.
  \end{equation}
  \begin{figure}[htbp]
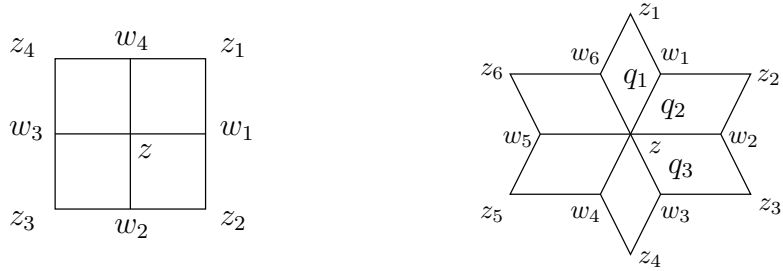

    \begin{center}
      \raisebox{0.4cm}{\input{squareLabels.tex}}\hfil{}
      \input{hexLabels.tex}
      \caption{Nine points of ${\mathbb Z}^2$ and 13 points of $\VTL$.}
      \label{fig:a:crLabeling}
    \end{center}
  \end{figure}
  To distinguish between the two sublattices, we will denote the
  points belonging to ${\mathbb Z}^2_0$ with $z_i$ while $w_i$ Are the
  points associated with ${\mathbb Z}^2_1$.
  Now (\ref{eq:a:crAForm}) gives four equations for a point $z$ and
  its eight neighbors as shown in Fig~\ref{fig:a:crLabeling}:
  \begin{equation}
    \label{eq:4crEquations}
    \begin{array}{rcl}\displaystyle
        (q - 1) \frac1{z -z_1}      &=& \displaystyle q\frac1{z - w_4} -
      \frac1{z-w_1},\\[0.4cm]
        \displaystyle(q^{-1} - 1) \frac1{z -z_2} &=&\displaystyle q^{-1}\frac1{z - w_1} -
      \frac1{z-w_2},\\[0.4cm]
        \displaystyle(q - 1) \frac1{z -z_3}      &=&\displaystyle q\frac1{z - w_2} -
      \frac1{z-w_3},\\[0.4cm]
        \displaystyle(q^{-1} - 1) \frac1{z -z_4} &=&\displaystyle q^{-1}\frac1{z - w_3} -
      \frac1{z-w_4}.
    \end{array}
  \end{equation}
  Multiplying the second and fourth equation with $q$ and taking the
  sum over all four equations (\ref{eq:4crEquations}) yields
  equation~(\ref{Toda(mn)}). Conversely given $w_1$, the first three
  equations (\ref{eq:4crEquations}) determine $w_2,w_3$, and $w_4$
  uniquely. The so determined $w_3$, and $w_4$ satify the fourth
  equation (\ref{eq:4crEquations}) (for any choice of $w_1$) if and
  only if (\ref{Toda(mn)}) holds.
  This proves the second part of the theorem.  Obviously interchanging
  $w$ and $z$ implies the claim for the other sublattice.
\end{proof}
Now we pass to the dual Kagome lattice shown in
Fig.~\ref{fig:QLattice}, where a similar relation holds for the
discrete equation of Toda type on the hexagonal lattice. However,
the symmetry between the sublattices is lost in this case.

Here we decompose $\VTL$ into $\VHL$ and
$\VTL\setminus\VHL\cong\FHL$.
\begin{theorem}
  Given $q_1,q_2,q_3\in {\mathbb C}$ with $q_1q_2q_3 = 1$ and a
  solution $z: \VTL\to \hat{\mathbb C}$ of the cross-ratio equations
  \begin{equation}
    \label{eq:a:cross-ratios}
    \begin{array}{rcl}
      q(z_{k,\ell,m-1},z_{k,\ell,m},z_{k,\ell+1,m},z_{k,\ell+1,m-1}) &=& q_1\\
      q(z_{k+1,\ell,m},z_{k,\ell,m},z_{k,\ell,m-1},z_{k+1,\ell,m-1}) &=& q_2\\
      q(z_{k,\ell-1,m},z_{k,\ell,m},z_{k+1,\ell,m},z_{k+1,\ell-1,m}) &=& q_3\\
    \end{array}
  \end{equation}
  then
  \begin{enumerate}
  \item restricted to $\FHL$ the solution $z$ satisfies the
    discrete equation of Toda type (\ref{eq.relToda}) on the
    hexagonal lattice
    \begin{equation}
      \label{eq:a:relToda}
      \sum_{k = 1}^3 A_k \left( \frac1{z - z_k} + \frac1{z -
          z_{k+3}}\right) = 0
    \end{equation}
  \item restricted to $\VHL$ the solution $z$ satisfies
    \begin{equation}
      \label{eq:a:crToda}
      \sum_{k = 1}^3 A_k \frac1{z-z_k} = 0
    \end{equation}
    where $A_i = \Delta_{i+2} -\Delta_{i+1}$ with $\Delta_i$ defined
    through $q_i = \frac{\Delta_{i+2}}{\Delta_{i+1}}$.
  \end{enumerate}
  Conversely given $q_1,q_2,q_3\in {\mathbb C}$ with $q_1q_2q_3 =
  1$, a solution $z$ to (\ref{eq:a:relToda}) on $\FHL$ and $z_{0,0,0}$,
  there is a unique extension $z:\VTL\to\hat{\mathbb C}$ satisfying
  (\ref{eq:a:cross-ratios}).

  Given $q_1,q_2,q_3\in {\mathbb C}$ with $q_1q_2q_3 =
    1$, a solution $z$ to (\ref{eq:a:crToda}) on $\VHL$ and
    $z_{1,0,0}$, there is a unique extension $z:\VTL\to\hat{\mathbb C}$
    satisfying (\ref{eq:a:cross-ratios}).

\end{theorem}
\begin{proof}
  First (\ref{eq:a:crAForm}) shows immediately that
  (\ref{eq:a:crToda}) is equivalent to having constant $\S$ as
  described in Section~\ref{sec:conformal}.

  Again we will distinguish the sublattices notationally by denoting
  points associated with elements of $\VHL$ with $w_i$ and points from
  $\FHL$ with $z_i$.  If the cross-ratios and neighboring points of a
  $z\in\FHL$ are labelled as shown in Fig.~\ref{fig:a:crLabeling},
  equation (\ref{eq:a:crAForm}) gives 6 equations:
  \begin{equation}
    \label{eq:6crEquations}
    \begin{array}{rcl}
    \displaystyle (q_1 - 1) \frac1{z -z_1} &=&\displaystyle q_1\frac1{z - w_6} - \frac1{z-w_1}\\
    \displaystyle (q_1 - 1) \frac1{z -z_4} &=&\displaystyle q_1\frac1{z - w_3} - \frac1{z-w_4}\\
    \displaystyle (q_2 - 1) \frac1{z -z_2} &=&\displaystyle q_2\frac1{z - w_1} - \frac1{z-w_2}\\
    \displaystyle (q_2 - 1) \frac1{z -z_5} &=&\displaystyle q_2\frac1{z - w_4} - \frac1{z-w_5}\\
    \displaystyle (q_3 - 1) \frac1{z -z_3} &=&\displaystyle q_3\frac1{z - w_2} - \frac1{z-w_3}\\
    \displaystyle (q_3 - 1) \frac1{z -z_6} &=&\displaystyle q_3\frac1{z - w_5} - \frac1{z-w_6}
    \end{array}
  \end{equation}
  To proof the first statement we take a linear combination of the
  equations (\ref{eq:6crEquations}). Namely $a$ the first two plus $b$
  times the second two plus $c$ times the third two. It is easy to see
  that there is a choice for $a,b$, and $c$ that makes the right hand
  side vanish if and only if $q_1 q_2 q_3 = 1$. If $q_1 =
  \frac{\Delta_3}{\Delta_2}, q_2 = \frac{\Delta_1}{\Delta_3},$ and
  $q_3 = \frac{\Delta_2}{\Delta_1}$ choose $a = \Delta_2, b=\Delta_3$,
  and $c = \Delta_1$. The remaining equation is (\ref{eq:a:relToda}).

  To proof the third statement we note that given $w_1$, we can
  compute $w_2$ through $w_6$ from equations (\ref{eq:6crEquations})
  2. to 6. but the ``closing condition'' equation
  (\ref{eq:6crEquations}) 1. is then equivalent to
  (\ref{eq:a:relToda}).

  The proofs of the second and fourth statements are literally the same
  if we choose $z_i = z_{i+3}$ and $w_i = w_{i+3}$.
\end{proof}

\section{Appendix. Proofs of Theorems of Section \ref{s.z^c}}\label{Appendix.lambda}

\subsection*{Proof of Theorem \ref{mu-theorem}}
Let $\Phi_{k,\ell,m}(\mu)$ be a solution to
(\ref{cross-ratio_delta/delta},\ref{constraint}) with some
$\mu$-independent matrices $C_{k,\ell,m}, B^{(n)}_{k,\ell,m}$. The
determinant identity
$$
\det \Phi_{k,\ell,m}(\mu)=(1-\mu \Delta_1)^k (1-\mu \Delta_2)^\ell
(1-\mu \Delta_3)^m \det \Phi_{0,0,0}(\mu)
$$
implies
$$
{\rm tr}\ \cA_{k,\ell,m}(\mu)= \dfrac{k}{\mu-1/\Delta_1}+
\dfrac{\ell}{\mu-1/\Delta_2}+ \dfrac{m}{\mu-1/\Delta_3}+a(\mu)
$$
with $a(\mu)$ independent of $k,\ell,m$. Without loss of
generality one can assume $a(\mu)=0$, i.e.
$$
{\rm tr}\ B^{(1)}_{k,\ell,m}=k, {\rm tr}\ B^{(2)}_{k,\ell,m}=\ell,
{\rm tr}\ B^{(3)}_{k,\ell,m}=m.
$$
This can be achieved by the change $\Phi\mapsto \exp(-1/2\int
a(\mu)d\mu)\Phi$.

The compatibility conditions of (\ref{wave evolution in mu}) and
(\ref{eq in mu}) read
\begin{eqnarray}                                        \label{compatibility}
\dfrac{d}{d\mu}\cL^{(1)}&=&\cA_{k+1,\ell,m}\cL^{(1)}-\cL^{(1)}\cA_{k,\ell,m}\nonumber\\
\dfrac{d}{d\mu}\cL^{(2)}&=&\cA_{k,\ell+1,m}\cL^{(2)}-\cL^{(2)}\cA_{k,\ell,m}\\
\dfrac{d}{d\mu}\cL^{(3)}&=&\cA_{k,\ell,m+1}\cL^{(3)}-\cL^{(3)}\cA_{k,\ell,m}\nonumber
\end{eqnarray}
The principal parts of these equations in $\mu=1/\Delta_1$ imply
\begin{eqnarray*}
B^{(1)}_{k+1,\ell,m} \left(\begin{array}{cc} 1
& f_1\\
\dfrac{1}{f_1} & 1\end{array}\right) &=& \left(\begin{array}{cc} 1
& f_1\\
\dfrac{1}{f_1} & 1\end{array}\right)
B^{(1)}_{k,\ell,m},\quad f_1=z_{k+1,\ell,m}-z_{k,\ell,m},\\
B^{(1)}_{k,\ell+1,m} \left(\begin{array}{cc} 1
& f_2\\
\dfrac{\Delta_2}{f_2\Delta_1} & 1\end{array}\right) &=&
\left(\begin{array}{cc} 1
& f_2\\
\dfrac{\Delta_2}{f_2\Delta_1 } & 1\end{array}\right)
B^{(1)}_{k,\ell,m},\quad f_2=z_{k,\ell+1,m}-z_{k,\ell,m},\\
B^{(1)}_{k,\ell,m+1} \left(\begin{array}{cc} 1
& f_3\\
\dfrac{\Delta_3}{f_3\Delta_1 } & 1\end{array}\right) &=&
\left(\begin{array}{cc} 1
& f_3\\
\dfrac{\Delta_3}{f_3\Delta_1} & 1\end{array}\right)
B^{(1)}_{k,\ell,m},\quad f_3=z_{k,\ell,m+1}-z_{k,\ell,m}.
\end{eqnarray*}
The solution with ${\rm tr}\ B^{(1)}_{k,\ell,m}=k$ is
$$
B^{(1)}_{k,\ell,m} =\dfrac{k-a_1}{z_{k+1,\ell,m}-z_{k-1,\ell,m}}
\left(\begin{array}{cc} z_{k+1,\ell,m}-z_{k,\ell,m}
& (z_{k+1,\ell,m}-z_{k,\ell,m})(z_{k,\ell,m}-z_{k-1,\ell,m})\\
1 & z_{k,\ell,m}-z_{k-1,\ell,m}\end{array}\right)+
\dfrac{a_1}{2}I.
$$
The same computation yields the formulas of Theorem
\ref{mu-theorem} for $B^{(2)}_{k,\ell,m}$ and
$B^{(3)}_{k,\ell,m}$.

To derive a formula for the coefficient $C_{k,\ell,m}$ let us
compare $\Phi_{k,\ell,m}(\mu)$ with  the solution
$\Psi_{k,\ell,m}(\lambda)$ in Theorem \ref{t.sym}, more exactly
with its extension to the lattice ${\Bbb Z}^3$:
\begin{eqnarray*}
\Psi_{k+1,\ell,m}(\lambda)&=&L^{(1)}(\lambda)\Psi_{k,\ell,m},\\
\Psi_{k,\ell+1,m}(\lambda)&=&L^{(2)}(\lambda)\Psi_{k,\ell,m},\\
\Psi_{k,\ell,m+1}(\lambda)&=&L^{(3)}(\lambda)\Psi_{k,\ell,m},
\end{eqnarray*}
normalized by $\Psi_{0,0,0}(\lambda)=I$. Here the matrices
$L^{(n)}$ are given by (\ref{L})
\begin{equation*}
L^{(n)}(\lambda)=(1-\lambda^2
\Delta_n)^{-1/2}\left(\begin{array}{cc}
1 & \lambda f_n\\
\lambda \Delta_n/f_n & 1 \end{array}\right).
\end{equation*}
Consider
\begin{equation*}
\tilde{\Psi}=h(\lambda) \left(\begin{array}{cc}
1/\sqrt\lambda & 0\\
0 & \sqrt\lambda \end{array}\right) \Psi \left(\begin{array}{cc}
\sqrt\lambda & 0\\
0 & 1/\sqrt\lambda \end{array}\right)
\end{equation*}
with
$$
h(\lambda)=(1-\lambda^2\Delta_1)^{k/2}(1-\lambda^2\Delta_2)^{\ell/2}(1-\lambda^2\Delta_3)^{m/2}.
$$
So defined $\tilde\Psi$ is a function of $\mu$. Since it satisfies
the same difference equations (\ref{wave evolution in mu}) as
$\Phi(\mu)$ and is normalized by
$$
\tilde\Psi_{k,\ell,m}(\mu=0)=I,
$$
we have
\begin{equation}                                \label{Phi-Psi}
\Phi_{k,\ell,m}(\mu)=\tilde\Psi_{k,\ell,m}(\mu)\Phi_{0,0,0}(\mu).
\end{equation}
Moreover $\tilde\Psi(\mu)$ is holomorphic in $\mu=0$ and due to
Theorem \ref{t.sym} equal in this point to
\begin{eqnarray*}
\tilde\Psi_{k,\ell,m}(\mu=0)=\left(\begin{array}{cc}
1 & Z\\
0 & 1\end{array}\right),\qquad Z=z_{k,\ell,m}-z_{0,0,0}.
\end{eqnarray*}
Taking the logarithmic derivative of (\ref{Phi-Psi}) with respect
to $\mu$
$$
\cA_{k,\ell,m}=\dfrac{d\tilde\Psi_{k,\ell,m}}{d\mu}\tilde\Psi^{-1}_{k,\ell,m}+
\tilde\Psi_{k,\ell,m} \cA_{0,0,0} \tilde\Psi^{-1}_{k,\ell,m}
$$
and computing its singularity at $\mu=0$ we get
\begin{eqnarray*}
C_{k,\ell,m}=\left(\begin{array}{cc}
1 & Z\\
0 & 1\end{array}\right) C_{0,0,0} \left(\begin{array}{cc}
1 & -Z\\
0 & 1\end{array}\right).
\end{eqnarray*}
The formula for $C_{k,\ell,m}$ in Theorem \ref{mu-theorem} is the
general solution to this equation.

Conversely, by direct computation one can check that the
compatibility conditions (\ref{compatibility}) with
$\cA_{k,\ell,m}$ computed above are equivalent to
(\ref{constraint}).

\subsection*{Proof of Theorem~\ref{t.compatibility}}

\begin{proof}
  First note that one can assume $b=0$ by applying a suitable
  M\"obius transformation. Next by translating $z$ one can make
  $z_{k,l,m} = 0$ for arbitrary fixed $k,l,m\in{\mathbb Z}$ (this will change
  $d$ however). Finally we can assume $d=0$ or $d=1$ since we
  can scale $z$.  Now one can show the compatibility by using a
  computer algebra system like MATHEMATICA as follows: given the
  points $z_{k,l,m}, z_{k+1,l,m}, z_{k,l\pm 1,m}, z_{k,l,m\pm 1}$ one
  can compute $z_{k-1,l,m}$ using the constraint (\ref{constraint}).
  With the cross-ratio equations one can now calculate all points
  necessary to apply the constraint for calculating $z_{k+2,l,m}$ and
  $z_{k,l,m+2}$ (besides $z_{k-1,l,m-1}$ one will need all points
  $z_{k\pm 1,l\pm 1,m\pm 1}$). Once again using the cross-ratios one
  calculates all points necessary to calculate $z_{k+1,l,m+2}$ with
  the constraint. Finally one can check, that the cross-ratio
  $q(z_{k,l,m+1},z_{k+1,l,m+1}, z_{k+1,l,m+2},z_{k,l,m+2})$ is
  correct. Since the three directions are equivalent and since all
  initial data was arbitrary (i.e. symbolic) this suffices to show the
  compatibility.
\end{proof}

\end{appendix}

\end{document}